\input amstex.tex
\magnification=1200
\documentstyle{amsppt}
\vcorrection{-1.0cm}
\pagewidth{32.5pc}
\topmatter
\title 
     Analysis on the minimal representation of $O(p,q)$
\\
     --  II. Branching laws
\endtitle
\author 
  Toshiyuki KOBAYASHI
  and 
  Bent \O RSTED
\endauthor
\affil 
RIMS Kyoto and SDU-Odense University
\endaffil
\address
RIMS, Kyoto University,
Sakyo-ku, Kyoto,  606-8502, Japan
\newline\indent
  Department of Mathematics and Computer Science,
  SDU-Odense University, Campusvej 55, DK-5230, Odense M, Denmark
\endaddress
\email{toshi\@kurims.kyoto-u.ac.jp;
       orsted\@imada.ou.dk}\endemail
\abstract
This is a second paper in a series devoted to the minimal
unitary representation of $O(p,q)$.  
By explicit methods from conformal geometry
 of pseudo Riemannian manifolds, 
 we find the branching law
corresponding to restricting the minimal unitary representation
to natural symmetric subgroups. In the case of
purely discrete spectrum we obtain the full spectrum
and give an explicit
Parseval-Plancherel formula, and in the
general case we construct an infinite discrete spectrum.

\endabstract
\endtopmatter
\NoRunningHeads
\overfullrule=0pt
\def\adm#1{{$#1$-admissible}}
\def\app{{+ +}}
\def\apm{{+ -}}

\def \Hsum#1{{{\underset{#1}\to{{\sum}^{\oplus}}}}}
\def \Ker{\operatorname{Ker}}

\def \trans{{}^t \!}
\def \tilLap#1{{\widetilde {\Delta}}_{#1}}
\def\V#1#2#3#4{V_{{#1}, {#4}}^{({#2}, {#3})}}

\def\Ad{\operatorname{Ad}}

\define \F#1#2{F({#1},{#2})}               %
\define \sgn{\operatorname{sgn}}

\def\ncone#1{{{\Cal N}_{#1}^*}}
\def \pro{{\operatorname{pr}}}
\def \Ass#1#2{{\Cal V_{{\frak {#1}}}({#2})}}

\def \Vpq{V^{p,q}}
\def\SSK#1{\operatorname{Supp}_K({#1})} %
\def\Kasym#1{\operatorname{AS}_K({#1})}
\def\pip#1#2#3{\pi^{{#1},{#2}}_{+, {#3}}}
\def\pim#1#2#3{\pi^{{#1},{#2}}_{-, {#3}}}
\def\der#1{\frac{\partial}{\partial{#1}}}
\define \set#1#2{\{{#1}:{#2}\}}
\def \sec#1{{\vskip 0.5pc\noindent$\underline{{\bold{\ch.\sc.}} \; \text{{#1}}}$\enspace}}
\def \num{{\ch.\sc}} %
\def\rarrowsim{\smash{\mathop{\,\rightarrow\,}\limits
  ^{\lower1.5pt\hbox{$\scriptstyle\sim$}}}}
\def\Pmax{P^{\roman{max}}}

\def\princeK#1#2{{\operatorname{Ind}_{\Pmax}^G({#2} \otimes \Bbb C_{#1})}}

\def\spr#1#2{{\Cal H^{#1}(\Bbb R^{#2})}}
\def \Dom{\operatorname{Dom}}
\def\kR#1#2#3{{\Xi(O({#1}) \times O({#2}):{#3})}}
\def\kK#1{{\Xi(K:{#1})}}
\define \zdfc_#1^#2#3{{\Cal R}_{\frak#1}^{#2} ({\Bbb C}_{#3})}
\define \zdf_#1^#2{{\Cal R}_{\frak {#1}}^{#2}}   %
\define \hc{{\frak h}^c}
\define \tc{{\frak t}^c}
\define \rootsys#1#2{\triangle(#1,#2)}
\define \metalk{(${\frak l}, (L\cap K)^\sim$)}
\define \Ato#1#2{{\Cal A}\left({#1}\triangleright{#2}\right)}

\define \possys#1#2{\triangle^+(#1,#2)}
\def \Adisc{{A'}}

\def \gk{$(\frak g, K)$}
\def \xbec{\cite{1}}
\def \xbz{\cite{2}}
\def \xbrkosI{\cite{3}}
\def \xbrkosII{\cite{4}}
\def \xbrkosIII{\cite{5}}
\def \xerdHigI{\cite{6}}
\def \xerdIntII{\cite{7}}
\def \xgs{\cite{8}}
\def \xhela{\cite{9}} %
\def \xhowe{\cite{10}}
\def \xhowetan{\cite{11}}
\def \xhuzhu{\cite{12}}
\def \xkobast{\cite{13}}
\def \xkupq{\cite{14}}
\def \xkdecomp{\cite{15}} %
\def \xkdecoalg{\cite{16}} %
\def \xkdecoass{\cite{17}} %
\def \xkmfjp{\cite{18}}
\def \xkmf{\cite{19}}
\def \xkdecoaspm{\cite{20}}
\def \xkhcrrest{\cite{21}}
\def \xkohcrcras{\cite{22}}
\def \xkorsI{\cite{23}}
\def \xkos{\cite{24}}
\def \xlee{\cite{25}}
\def \xsab{\cite{26}}
\def \xschlap{\cite{27}}
\def \xschm{\cite{28}}%
\def \xschmid{\cite{}}%
\def \xtoramin{\cite{29}}
\def        \xvg{\cite{30}}%
\def \xvu{\cite{31}}
\def        \xvr{\cite{32}}%
\def \xvi{\cite{33}}
\def      \xwong{\cite{34}}%
\def \xorslmp{\cite{35}}
\def \xorsjfa{\cite{36}}

\document

\centerline{{\bf {Contents}}}
\smallskip
\item{}\enspace Introduction
\item{\S 4.}\enspace
 Criterion for discrete decomposable branching laws
\item{\S 5.}\enspace
 Minimal elliptic representations of $O(p,q)$ 
\item{\S 6.}\enspace
  Conformal embedding of the hyperboloid
\item{\S 7.}\enspace
  Explicit branching formula (discretely decomposable case)
\item{\S 8.}\enspace 
 Inner product on $\varpi^{p,q}$ and the Parseval-Plancherel formula
\item{\S 9.}\enspace
  Construction of discrete spectra in the branching laws

\head
Introduction
\endhead

This is the second in a series of papers devoted to the
analysis of the minimal representation $\varpi^{p,q}$ of $O(p,q)$. We refer
to \xkorsI{}  for a general introduction; also the
numbering of the sections is continued from that paper,
and we shall refer back to sections there. 
However, 
 the present paper may be read independently from \xkorsI, and
 and its main object is to study the branching law for the minimal
 unitary representation $\varpi^{p,q}$
 from {\bf analytic} and {\bf geometric} point of view. 
 Namely, we shall find
 by explicit means, coming from conformal geometry,
the restriction of $\varpi^{p,q}$ with respect to the symmetric pair
$$
  (G, G') = (O(p,q), O(p',q')\times O(p'',q'')).
$$
 If one of
the factors in $G'$ is compact, then the spectrum is discrete 
 (see Theorem~4.2 also for an opposite implication), and
we find the explicit branching law; when both factors are non-compact,
there will still (generically) be an infinite discrete spectrum, which we
also construct (conjecturally almost all of it; see \S 9.8).
We shall see that the (algebraic)
situation is similar to the theta-correspondence, where the
metaplectic representation is restricted to analogous subgroups.

Let us here state the main results in a little more precise form,
referring to sections 8 and 9 for further notation and details.

\proclaim{Theorem~A \
{\rm (the branching law for
 $O(p,q) \downarrow O(p,q') \times O(q'')$; see Theorem~7.1)}}
If $q'' \ge 1$ and $q' + q'' = q$,
 then
the twisted pull-back $\widetilde{\Phi_1^*}$
 of the local conformal map $\Phi_1$
 between spheres
 and hyperboloids
 gives an explicit irreducible decomposition
 of the unitary representation $\varpi^{p,q}$
 when restricted to $O(p,q') \times O(q'')$:
$$
  \widetilde{{(\Phi_1})^*}
 \: 
 \varpi^{p,q}|_{O(p, q') \times O(q'')}
 \rarrowsim
 \overset{\infty\hphantom{M}}\to{\Hsum{l=0\hphantom{M}}}
  \pip{p}{q'}{l + \frac{q''}2 -1} \boxtimes \spr{l}{q''}.
$$
\endproclaim

 The representations appearing in the decompositions
are in addition to usual spherical harmonics $\spr{l}{q}$ for compact
orthogonal groups $O(q)$, also the representations $\pip{p}{q}{\lambda}$ 
for non-compact orthogonal groups $O(p,q)$. 
The latter ones may be thought of 
as  discrete series representations on hyperboloids 
$$
   X(p,q) :=\set{x=(x',x'') \in \Bbb R^{p+q}}{|x'|^2-|x''|^2=1}
$$
for $\lambda >0$
 or their analytic continuation for $\lambda \le 0$;
 they may be also thought of
 as cohomologically induced representations from
 characters of  certain $\theta$-stable
  parabolic subalgebras. 
The fact that they occur in this branching law gives
a different proof of the unitarizability
 of these modules $\pip{p}{q}{\lambda}$
 for $\lambda > -1$,
 once we know $\varpi^{p,q}$ is unitarizable (cf. Part I, Theorem~3.6.1).
It might be interesting to remark that the unitarizability 
  for $\lambda < 0$
 (especially, $\lambda =  -\frac12$ in our setting) 
 does follow neither from a general unitarizability theorem
  on Zuckerman-Vogan's derived functor modules \xvu,
 nor from a general theory of harmonic analysis on semisimple symmetric spaces.
 
Our intertwining operator $\widetilde{{(\Phi_1})^*}$ in Theorem~A is 
 derived from a conformal change of coordinate (see \S 6 for its explanation)
  and is explicitly written.
 Therefore,
 it makes sense to ask also about the relation of unitary inner products
  between the left-hand and the right-hand side in the branching formula.
Here is an answer (see Theorem~8.5):  
We normalize the inner product $\| \ \|_{\pip{p}{q}{\lambda}}$ (see \S 8.3)
 such that for $\lambda>0$,
$$
    \| f\|_{\pip{p}{q}{\lambda}}^2  = 
  \lambda \| f\|_{L^2(X(p,q))}^2,
  \quad
  \text{ for any } f \in (\pip{p}{q}{\lambda})_K.  
$$
\proclaim{Theorem~B \
{\rm (the Parseval-Plancherel formula for
 $O(p,q) \downarrow O(p,q') \times O(q'')$)}}
\newline
{\rm 1)}\enspace
If we develop $F \in \Ker{\tilLap{M}}$ as
 $F = \sum_l^\infty F_{l}^{(1)} F_{l}^{(2)}$
 according to the irreducible decomposition in Theorem~A,
 then we have
$$
  \|F\|_{\varpi^{p,q}}^2
 = \sum_{l=0}^\infty
  \| F_{l}^{(1)}\|_{\pip{p}{q'}{l+\frac{q''}2-1}}^2 \|
  \
  F_{l}^{(2)}\|_{L^2(S^{q''-1})}^2.
$$
\noindent
{\rm 2)}\enspace
In particular,
 if $q'' \ge 3$,
 then all of $\pip{p}{q'}{l + \frac{q''}2-1}$ are
 discrete series for the hyperboloid $X(p,q')$
 and the above formula amounts to
$$
  \|F\|_{\varpi^{p,q}}^2
 = \sum_{l=0}^\infty
  (l + \frac{q''}2-1)
  \
  \| F_{l}^{(1)}\|_{L^2(X(p,q'))}^2 \|
  \
  F_{l}^{(2)}\|_{L^2(S^{q''-1})}^2.
$$
\endproclaim

The formula may be also regarded as an explicit unitarization
 of the minimal representation $\varpi^{p,q}$ on the \lq\lq hyperbolic space model\rq\rq\
 by means of the right side
 (for an abstract unitarization of $\varpi^{p,q}$,
  it suffices to choose a single
 pair $(q',q'')$).
We note that the formula was previously known
 in the case where $(q', q'') = (0,q)$
 (namely, when each summand in the right side is finite dimensional)
 by Kostant, Binegar-Zierau by a different approach.
The formula is new and seems to be particularly interesting even in the
 special case $q''=1$,
 where the minimal representation $\varpi^{p,q}$
 splits into two irreducible summands
 when restricted to $O(p,q-1)\times O(1)$.

In Theorem~9.1, we consider a more general setting and prove:
\proclaim{Theorem~C\
{\rm (discrete spectrum in the restriction
 $O(p,q) \downarrow O(p' ,q') \times O(p'',q'')$)}}
\newline
The twisted pull-back of the locally conformal diffeomorphism
 also constructs
$$
\Hsum{\lambda \in \Adisc(p', q') \cap \Adisc(q'',p'')}
  \pip{p'}{q'}{\lambda} \boxtimes \pim{p''}{q''}{\lambda}
\oplus
\Hsum{\lambda \in \Adisc(q', p') \cap \Adisc(p'',q'')}
  \pim{p'}{q'}{\lambda} \boxtimes \pip{p''}{q''}{\lambda}
$$
as a discrete spectrum in the branching law for the non-compact case.
\endproclaim

Even in the special case $(p'',q'')=(0,1)$,
 our branching formula includes a new and mysterious construction
  of the minimal representation on the hyperboloid as below (see Corollary~7.2.1):
Let $W^{p,r}$ be the set of $K$-finite vectors ($K = O(p) \times O(r)$) of
 the kernel of the Yamabe operator
$$
 \Ker \tilLap{X(p,r)}= \set{f \in C^\infty(X(p,r))}{\Delta_{X(p,r)} f = \frac14 {(p+r-1)(p+r-3)} f},
$$
 on which the isometry group $O(p,r)$ 
 and the Lie algebra of the conformal group $O(p,r+1)$ act.
 The following Proposition is a consequence
  of Theorem~7.2.2 by an elementary linear algebra.
\proclaim{Proposition~D}
Let $m > 3$ be odd.
There is a long exact sequence
$$
 0 \rightarrow W^{1,m-1} \overset{\varphi_1}\to\longrightarrow
 W^{2,m-2} \overset{\varphi_2}\to\longrightarrow
  W^{3,m-3} \overset{\varphi_3}\to\longrightarrow
 \dotsb
  \overset{\varphi_{m-2}}\to\longrightarrow
  W^{m-1,1} \overset{\varphi_{m-1}}\to\longrightarrow
  0
$$
such that
 $\Ker \varphi_{p}$ is isomorphic to $(\varpi^{p,q})_K$
  for any $(p,q)$ such that $p+q=m+1$.
\endproclaim
We note that each representation space $W^{p,q-1}$ 
 is realized on a different space
 $X(p,q-1)$ whose isometry group $O(p,q-1)$ varies
  according to $p$ ($1 \le p \le m$).
So, one may expect that only the intersections of adjacent groups can act
 (infinitesimally) on $\Ker \varphi_p$.
Nevertheless, a larger group $O(p,q)$ can act on 
 a suitable completion of $\Ker \varphi_p$, 
 giving rise to another construction of the minimal representation 
 on the hyperboloid $X(p,q-1) = O(p,q-1)/O(p-1,q-1)$ !
We note that
 $\Ker \varphi_p$ is roughly
  half the kernel of the Yamabe operator on the hyperboloid
  (see \S 7.2 for details).

We briefly indicate the contents of the paper:
In section 4 we recall the relevant facts about discretely decomposable
restrictions from  \xkdecoalg\ and \xkdecoass,
 and apply the criteria to our present situation. In
particular, we calculate the associated variety of $\varpi^{p,q}$
as well as its asymptotic $K$-support introduced by Kashiwara-Vergne.
Theorem~4.2 and Corollary~4.3 clarify
 the reason why we start with the subgroup $G'=O(p,q') \times O(q'')$
 (i.e. $p''=0$).
 Section 5 contains the identification
of the representations $\pip{p}{q}{\lambda}$ and $\pim{p}{q}{\lambda}$
of $O(p,q)$ in several ways, namely as: derived functor modules,
Dolbeault cohomologies, eigenspaces on hyperboloids, and quotients or 
subrepresentations of parabolically induced modules. In section 6 we give
the main construction of embedding conformally a direct product
of hyperboloids into a product of spheres; this gives rise to
a canonical intertwining operator between solutions to the so-called 
Yamabe equation, studied in connection with conformal differential geometry,
on conformally related spaces. Applying this principle in section 7 we
obtain the branching law in the case where one factor in $G'$ is compact,
and in particular when one factor is just $O(1)$. In this case we have
Corollary 7.2.1, stating that $\varpi^{p,q}$ restricted to $O(p,q-1)$  
is the direct sum of two representations, realized in even respectively odd
functions on the hyperboloid for $O(p,q-1)$. Note here the analogy
with the metaplectic representation. Also note here
Theorem 7.2.2, which gives a mysterious extention of $\varpi^{p,q}$
by $\varpi^{p+1,q-1}$ - both inside the space of solutions to
the Yamabe equation on the hyperboloid $X(p,q-1) = O(p,q-1)/O(p-1,q-1)$. 
We also point out that the representations $\pip{p}{q}{\lambda}$ for
$\lambda = 0, -\frac 1 2$ are rather exceptional; they are unitary,
but outside the usual \lq\lq fair range" for derived functor modules, see
the remarks in section 8.2.  
Section 8 contains a calculation
of spectra of intertwining operators and gives the explicit normalization
of a Parseval-Plancherel formula for the branching law, and finally
in section 9 we use certain Sobolev estimates to construct an infinite
discrete spectrum when both factors in $G'$ are non-compact. We also conjecture
the form of the full discrete spectrum (true in the case of a compact factor).
It should be interesting to calculate the full Parseval-Plancherel formula
in the case of both discrete and continuous spectrum. 

The first author expresses his sincere gratitude to SDU - Odense University for
the warm hospitality.

\def \ch{4}
\def \sc{1}
\head
\S \ch. Criterion for discrete decomposable branching laws
\endhead

\def \sc{1}
\sec{}
Our object of study is the discrete spectra
 of the branching law of the restriction $\varpi^{p,q}$
 with respect to a symmetric pair
 $(G, G') = (O(p,q), O(p',q')\times O(p'',q''))$.
The aim of this section is to
 give a necessary and sufficient condition on $p', q', p''$ and $q''$
 for the branching law to be discretely decomposable.

We start with general notation.
Let $G$ be a linear reductive Lie group,
 and $G'$ its subgroup which is reductive in $G$. 
We take a maximal compact subgroup $K$ of $G$
 such that $K' := K \cap G'$ is also a maximal compact subgroup.
Let $\frak g_0 = \frak k_0 + \frak p_0$ be a Cartan decomposition,
 and $\frak g = \frak k + \frak p$ its complexification.
Accordingly, we have a direct decomposition
 $\frak g^* = \frak k^* + \frak p^*$
 of the dual spaces.
 
Let $\pi \in \widehat{G}$.
We say that the restriction $\pi|_{G'}$ is {\it \adm{G'}}
 if $\pi|_{G'}$ splits into a direct Hilbert sum
 of irreducible unitary representations of $G'$
 with each multiplicity finite (see \xkdecomp).
As an algebraic analogue of this notion, 
 we say the underlying $(\frak g, K)$-module
 $\pi_K$ is {\it discretely decomposable as a $(\frak g', K')$-module},
 if $\pi_K$ is decomposed into
 an algebraic direct sum of irreducible $(\frak g', K')$-modules
 (see \xkdecoass).
We note that if the restriction $\pi|_{K'}$ is \adm{K'},
 then the restriction $\pi|_{G'}$ is also \adm{G'} (\xkdecomp, Theorem~1.2)
 and the underlying \gk-module $\pi_K$ is discretely decomposable
 (see \xkdecoass, Proposition~1.6).
Here are criteria for $K'$-admissibility and discrete decomosability:
\proclaim{Fact \num\ {\rm (see \xkdecoalg, Theorem~2.9 for (1);
 \xkdecoass, Corollary~3.4 for (2))}}
\newline
{\rm 1)} If $\Kasym{\pi} \cap \Ad^*(K) (\frak k')^{\perp} = 0$,
  then $\pi$ is \adm{K'} and also \adm{G'}.
\newline
{\rm 2)} If $\pi_K$ is discretely decomposable as a $(\frak g',K')$-module,
 then
$
 \pro_{\frak g \to \frak g'} \left(\Ass{g}{\pi_K}\right)
 \subset   \ncone{\frak g'}.
$
\endproclaim
Here,
 $\Kasym{\pi}$ is the asymptotic cone of
$$\SSK{\pi}:=
 \set{\text{highest weight of }\tau \in \widehat{K_0}}{
    [\pi|_{K_0}:\tau]\neq 0}
$$
 where $K_0$ is the identity component of $K$, 
 and
 $(\frak k')^{\perp} \subset \frak k^*$ is the annihilator of $\frak k'$.
Let $\ncone{\frak p} \ (\subset \frak p^*)$ be
 the nilpotent cone for $\frak p$.
$\Ass{g}{\pi_K}$ denotes the associated variety of $\pi_K$,
 which is an $\Ad^*(K_\Bbb C)$-invariant closed subset of $\ncone{\frak p}$.
We write the projection
 $\pro_{\frak p\to \frak p'}\: \frak p^* \to {\frak p'}^*$ 
 dual to the inclusion   $\frak p' \hookrightarrow \frak p$.

\def \sc{2}
\sec{}
Let us consider our setting where
 $\pi = \varpi^{p,q}$
 and $(G, G') = (O(p,q), O(p',q') \times O(p'', q''))$.
\proclaim{Theorem~\num}
Suppose $p' + p'' =p \ (\ge 2)$, $q' + q'' = q \ (\ge 2)$ and
 $p+q \in 2 \Bbb N$.
Then the following three conditions on $p', q', p'', q''$ are equivalent:
\newline {\rm i)} \enspace
$\varpi^{p,q}$ is \adm{K'}.
\newline {\rm ii)} \enspace
$\varpi^{p,q}_K$ is discretely decomposable as a $(\frak g', K')$-module.
\newline {\rm iii)} \enspace
$\min(p', q', p'', q'') = 0$.
\endproclaim

The implication (i) $\Rightarrow$ (ii) holds by a general
 theory  as we explained (\xkdecoass,  Proposition~1.6);
 (ii) $\Rightarrow$ (iii) will be proved in \S \ch.4,
 and
 (iii) $\Rightarrow$ (i) in \S \ch.5,
 by an explicit computation of the asymptotic cone
 $\Kasym{\varpi^{p,q}}$ and the associated variety
$\Ass{g}{\varpi^{p,q}_K}$ which are used in Fact~\ch.1.

\remark{Remark} %
Analogous results to the equivalence (i) $\Leftrightarrow$ (ii) in Theorem~\num\
 were first proved in \xkdecoass, Theorem~4.2
 in the setting where $(G, G')$ is any reductive symmetric pair and the representation
  is any $A_{\frak q}(\lambda)$ module in the sense of Zuckerman-Vogan,
  which may be regarded as \lq\lq representations attached to elliptic orbits\rq\rq .
We note that our representations $\varpi^{p,q}$ are supposed to be attached to nilpotent orbits.
We refer \xkdecoaspm, Conjecture~A to relevant topics.
\endremark
\def \sc{3}
\sec{}
The following corollary is a direct consequence of Theorem~\ch.2,
 which will be an algebraic background
 for the proof of the explicit branching law
 (Theorem~7.1).
\proclaim{Corollary~\num}
Suppose that  $\min(p', q', p'', q'') =0$.
\newline
{\rm 1)}\enspace 
The restriction of the unitary representation
 $\varpi^{p,q}|_{G'}$ is also \adm{G'}.
\newline
{\rm 2)}\enspace 
 The space of $K'$-finite vectors $\varpi^{p,q}_{K'}$ coincides
 with that of $K$-finite vectors $\varpi^{p,q}_K$.
\endproclaim
\demo{Proof} See \xkdecomp, Theorem~1.2 for (1);
 and \xkdecoass, Proposition~1.6 for (2).
\qed
\enddemo

A geometric counterpart of Corollary~\num~(2) is reflected as
 the removal of singularities of matrix coefficients
 for the discrete spectra in the analysis that we study in \S 6;
 namely,
 any analytic function defined on an open subset $M_+$ 
 (see \S 6 for notation) of $M$
 which is a $K'$-finite vector of a discrete spectrum,
 extends analytically on $M$ 
 if $p'' = 0$.
The reason for this is not only the decay of matrix coefficients
 but a matching condition of the leading terms for $t \to \pm \infty$.
This is not the case for $\min(p', q', p'', q'') > 0$ (see \S 9).

\def \sc{4}
\sec{}
Proof of (ii) $\Rightarrow$ (iii) in Theorem~\ch.2.

We identify $\frak p^*$ with $\frak p$ via the Killing form,
 which is in turn identified with $M(p,q;\Bbb C)$ by
$$
   M(p,q;\Bbb C) \rarrowsim \frak p,
  \
   X \mapsto \pmatrix O & X \\ \trans{X} & O \endpmatrix.
$$
Then the nilpotent cone $\ncone{\frak p}$ corresponds to the following variety.
$$
  \set{X \in M(p,q;\Bbb C)}{\text{ both }
     X \trans{X} \text{ and }  \trans{X} X
    \text{ are nilpotent matrices}}
\tag \num.1
$$
We put
$$
 M_{0,0}(p,q;\Bbb C) := 
  \set{X \in M(p,q;\Bbb C)}{X \trans{X} = O, \trans{X} X = O}.
$$
Then $M_{0,0}(p,q;\Bbb C)\setminus\{O\}$ 
 is the unique
 $K_\Bbb C \simeq O(p, \Bbb C) \times O(q, \Bbb C)$-orbit of dimension $p+q-3$.
The associated variety $\Ass{g}{\varpi^{p,q}_K}$
 of $\varpi^{p,q}$ is of dimension $p+q-3$,
 which follows easily from the $K$-type formula of $\varpi^{p,q}$
 (see \xkorsI, Theorem~3.6.1).
Thus,  we have proved:
\proclaim{Lemma~\num}
The associated variety $\Ass{g}{\varpi^{p,q}_K}$
 equals $M_{0,0}(p,q;\Bbb C)$.
\endproclaim

The projection $\pro_{\frak p \to \frak p'} \: \frak p^* \to {\frak p'}^*$ 
 is identified with the map
$$
\pro_{\frak p \to \frak p'}\:
  M(p,q;\Bbb C) \to M(p', q'; \Bbb C) \oplus M(p'', q''; \Bbb C),
   \
  \pmatrix X_1 & X_2 \\ X_3 & X_4 \endpmatrix
  \mapsto (X_1, X_4).
$$
Suppose $p' p'' q' q'' \neq 0$.
If we take 
$$
 X := E_{1, 1} - E_{p'+1, q'+1} 
 + \sqrt{-1} E_{p'+1, 1} + \sqrt{-1} E_{1, p'+1} \in M_{0,0}(p,q;\Bbb C),
$$
 then $\pro_{\frak p \to \frak p'} (X) = (E_{1, 1}, - E_{p'+1, q'+1})$.
But $E_{1,1} \not\in \ncone{\frak{o}(p',q')}$ and
  $-E_{p'+1,q'+1} \not\in \ncone{\frak{o}(p'',q'')}$.
Thus, 
 $\pro_{\frak g \to \frak g'} (X) \not\in   \ncone{\frak g'}$.
It follows from Fact~\ch.1~(2) that $\varpi^{p,q}_K$
 is not discrete decomposable as a $(\frak g', K')$-module.
Hence (ii) $\Rightarrow$ (iii) in Theorem~\ch.2 is proved.
\qed
\def \sc{5}
\sec{}
Proof of (iii) $\Rightarrow$ (i) in Theorem~\ch.2.

We take an orthogonal complementary subspace $\frak k_0''$ of $\frak k_0'$
 in $\frak k_0 \simeq {\frak o}(p) + \frak o(q)$.
Let $\tc_0$ be a Cartan subalgebra of $\frak k_0$
 such that $\frak t_0'' := \tc_0 \cap \frak k_0''$
 is a maximal abelian subspace in $\frak k_0''$.
We choose a positive system $\Delta^+(\frak k, \tc)$
 which is compatible with a positive system of
 the restricted root system $\Sigma(\frak k, \frak t'')$.
Then we can find a basis $\set{f_i}{1 \leq i \leq [\frac{p}2]+[\frac{q}2]}$
 on $\sqrt{-1}\frak t_0^*$
 such that a positive root system of $\frak k$ is given by
$$
\align
    \possys{\frak k}{\tc}
    =
    &\set{f_i \pm f_j}{1 \leq i < j \leq[\frac p2]}
\\
    &\cup
    \set{f_i \pm f_j}{[\frac {p}2] + 1 \leq i < j \leq [\frac{p}2]+[\frac{q}2]}
\\
    &\cup
    \left( \set{\pm f_l}{1 \leq l \leq [\frac p2]} 
    \;{\text { ($p$:odd) }}
    \right)
\\
    &\cup
    \left( \set{\pm f_l}{[\frac{p}2]+1 \leq l \leq [\frac{p}2]+[\frac{q}2] } 
    \;{\text { ($q$:odd) }}
    \right),
\endalign
$$
and such that
$$
 \sqrt{-1} (\frak t_0'')^* 
 = 
 \sum_{j=1}^{\min(p', p'')}   \Bbb R f_{i}
 +
 \sum_{j=1}^{\min(q', q'')}   \Bbb R f_{[\frac{p}{2}]+j}
\tag \num.1
$$
 if we regard $(\frak t_0'')^*$ as a subspace of $(\tc_0)^*$
 by the Killing form.

Suppose $p' q' p'' q'' = 0$.
Without loss of generality
 we may and do assume $p'' = 0$,
 namely,
 $G' = O(p,q') \times O(q'')$ with $q' + q'' = q$.

Let us first consider the case $p \neq 2$.
Then the irreducible $O(p)$-representation $\spr{a}{p}$
 remains irreducible when restricted to $SO(p)$.
The corresponding highest weight is given by $a f_1$.
It follows from the $K$-type formula of $\varpi^{p,q}$ 
 (Theorem~3.6.1)
 that
$$
\SSK{{\varpi}^{p,q}} = \set{a f_1 + b f_{[\frac{p}2]+1}}{a, b \in \Bbb N, 
 a + \frac{p}2 = b + \frac{q}2}.
$$
Therefore we have proved:
$$
    \Kasym{{\varpi}^{p,q}} = \Bbb R_+ (f_1 + f_{[\frac{p}{2}]+1}).
\tag \num.2
$$
Then
$\Kasym{{\varpi}^{p,q}} \cap \sqrt{-1}(\frak t_0'')^* = \{0\}$
 from (\num.1) and (\num.2),
 which implies
$$
 \Kasym{{\varpi}^{p,q}} \cap \sqrt{-1} \Ad^*(K)(\frak k'_0)^\perp = \{0\}
$$
 because $(\frak t_0'')^*$ meets any $\Ad^*(K)$-orbit through
 $(\frak k'_0)^\perp$.
Therefore,
 the restriction ${\varpi}^{p,q}|_{K'}$ is \adm{K'}
 by Fact~\ch.1~(1).

If $p=2$, then $\varpi^{p,q}$ splits into two representations
 (see Remark~3.7.3), 
 say $\varpi^{2,q}_+$ and $\varpi^{2,q}_-$,
 when restricted to the connected component $SO_0(2,q)$.
Likewise, $\spr{a}{p}$ is a direct sum of two one-dimensional representations
 when restricted to $SO(2)$ if $a \ge 1$.
Then we have 
$$
    \Kasym{{\varpi}^{2,q}_\pm} = \Bbb R_+ (\pm f_1 + f_{p+1}).
$$
Applying Fact~4.1~(1) to the identity components $(G_0, G'_0)$
 of groups $(G, G')$,
 we conclude that
 the restriction ${\varpi}_\pm^{p,q}|_{K'_0}$ is \adm{K'_0}.
Hence
 the restriction ${\varpi}^{p,q}|_{K'}$ is also \adm{K'}.
Thus, (iii) $\Rightarrow$ (i) in Theorem~\ch.2 is proved.

Now the proof of Theorem~\ch.2 is completed.
\qed

\def \ch{5}
\def \sc{1}
\head
 \S \ch.  Minimal elliptic representations of $O(p,q)$
\endhead

\def \sc{1}
\sec{}
In this section,
 we introduce a family of irreducible representations of $G = O(p,q)$,
 denoted by  $\pip{p}{q}\lambda, \pim{p}{q}\lambda$,
 for $\lambda \in A_0(p,q)$,
 in three different realizations.
These representations %
 are supposed to be attached to minimal elliptic orbits,
 for $\lambda > 0$
 in the sense of the Kirillov-Kostant orbit method.
Here, we set
$$
        A_0(p,q) := 
        {\cases
        {\set {\lambda \in \Bbb Z + \frac {p+q} 2} {\lambda >-1}}
        \quad
        &(p>1, q \ne 0),
\\
         {\set{\lambda \in \Bbb Z + \frac {p+q} 2} {\lambda \ge \frac p 2 -1}}
        &(p>1, q = 0),
\\
         \emptyset  
        &(p=1,q \ne 0) \text{ or } (p=0),
\\
         \{-\frac 1 2, \frac 1 2\} 
        &(p=1,q=0).
        \endcases}
\tag \num.1
$$
It seems natural to include the parameter $\lambda = 0, -\frac12$
 in the definition of $A_0(p,q)$ as above, 
 although $\lambda = -\frac12$ is outside the weakly fair range
 parameter in the sense of Vogan \xvi.
Cohomologically induced representations for $\lambda = -\frac12$
 and $\lambda = -1$ will be discussed in details in a subsequent paper.
In particular, the case $\lambda = -1$ is of importance
 in another geometric construction of the minimal representation
 via Dolbeault cohomology groups (see Part~I, Introduction, Theorem~B~(4)).

\def \sc{2}
\sec{}
Let $\Bbb R^{p,q}$ be the Euclidean space $\Bbb R^{p+q}$
 equipped with the flat pseudo-Riemannian metric:
$$
 g_{\Bbb R^{p,q}} = {d v_0}^2 + \cdots + d {v_{p-1}}^2 - {d v_{p}}^2 
         - \cdots - d {v_{n+1}}^2.  
$$
We define a hyperboloid by
$$
     X(p,q):= \set{(x,y) \in \Bbb R^{p,q}}
                  {|x|^2 - |y|^2=1}.
$$
We note $X(p,0) \simeq S^{p-1}$ and $X(0,q) = \emptyset$.  
If $p=1$,
 then $X(p,q)$ has two connected components.  
The group $G$ acts transitively on $X(p,q)$
 with isotropy subgroup $O(p-1,q)$ at 
$$
x^o := \trans (1,0, \cdots, 0).
\tag \num.1
$$  
Thus $X(p,q)$ is realized as a homogeneous manifold:
$$
     X(p,q) \simeq O(p,q) /O(p-1,q).  
$$
We induce a pseudo-Riemannian metric $g_{X(p,q)}$ on $X(p,q)$
 from $\Bbb R^{p,q}$ (see \xkorsI, \S 3.2),
 and write $\Delta_{X(p,q)}$
 for the Laplace-Beltrami operator on $X(p,q)$.  
As in \xkorsI, Example~2.2,
 the Yamabe operator is given by
$$
     \tilLap{X(p,q)}= \Delta_{X(p,q)} -\frac14 {(p+q-1)(p+q-3)}.
\tag \num.2
$$
For $\lambda \in \Bbb C$,
 we set
$$
\alignat1
   C_\lambda^\infty(X(p,q))
 &:= \set{f \in C^\infty(X(p,q))}{
       \Delta_{X(p,q)} f =(-\lambda^2 + \frac14 (p+q -2)^2)f}
\\
 & = \set{f \in C^\infty(X(p,q))}{
       \tilLap{X(p,q)} f =(-\lambda^2 + \frac 1 4)f}.
\tag \num.3
\endalignat
$$
Furthermore, 
 for $\epsilon = \pm $,
 we write
$$
   C_\lambda^\infty(X(p,q))_\epsilon
  := \set{f \in C_\lambda^\infty(X(p,q))}{
f(-z) = \epsilon f(z), \ z \in X(p,q)}.
$$
Then we have a direct sum decomposition
$$
 C_\lambda^\infty(X(p,q)) =
   C_\lambda^\infty(X(p,q))_+ + 
   C_\lambda^\infty(X(p,q))_-
\tag \num.4
$$
 and each space is invariant under left translations of
 the isometry group $G$
 because $G$ commutes with $\Delta_{X(p,q)}$.
With the notation in \S 3.5, 
 we note
 if  $q = 0$, then 
 $C^\infty_\lambda(X(p,0))_{\sgn (-1)^k}$ is finite dimensional
 and isomorphic to the space of spherical harmonics:
$$
  \spr{k}{p} \simeq C^\infty_\lambda(X(p,0))
 = C^\infty_\lambda(X(p,0))_{\sgn (-1)^k}
 \quad \text{ $(k := \lambda + \frac{p-2}2)$}.
$$

\def \sc{3}
\sec{}
Let $G = O(p,q)$ where $p, \,q \geq 1$
and let $\theta$ be the Cartan involution corresponding to 
$K = O(p) \times O(q)$.
We extend a Cartan subalgebra $\tc_0$ of $\frak k_0$ (given in \S 4.5)
 to that of $\frak g_0$, denoted by $\hc_0$.
If both $p$ and $q$ are odd, 
 then $\dim \hc_0 = \dim \tc_0 +1$;
 otherwise $\hc_0 = \tc_0$.
The complexification of $\hc_0$ is denoted by $\hc$.

We can take
 a basis $\set{f_i}{1 \leq i \leq [\frac {p+q}2]}$ of $(\hc)^*$
 (see \S 4.5; by a little abuse of notation if both $p$ and $q$ are odd)
 such that the root systems of $\frak g$ is given by
$$ 
\aligned
    \rootsys{\frak g}{\hc} 
    = 
    &\set{\pm (f_i \pm f_j)}{1 \leq i < j \leq [\frac{p+q}2]}
\\
    &\cup \left( \set{ \pm f_l}{1 \leq l \leq [\frac{p+q}2]}
    \;{\text { ($p+q$:odd) }} \right).
\endaligned
$$

Let $\{ H_i \} \subset \hc$ be the dual basis for $\{ f_i \} \subset (\hc)^*$. 
Set ${\frak t} :=  {\Bbb C} H_1$  ($\subset {\frak t}^c \subset \hc$).
Then  the centralizer $L$ of ${\frak t}$ in $G$
 is isomorphic to $SO(2) \times O(p-2,q)$.
Let ${\frak q} = {\frak l} + {\frak u}$
 be a $\theta$-stable parabolic subalgebra of $\frak g$
 with nilpotent radical ${\frak u}$ given by,
$$
    \rootsys {\frak u} \hc 
    := 
   \set{f_1 \pm f_j}{2 \le j \le [\frac{p+q}2]}
   \cup
   \left( \{ \pm f_1\}  \;{\text { ($p+q$:odd) }} \right),
$$
 and with a Levi part $\frak l = \frak l_0 \otimes \Bbb C$
 given by
$$
  \frak l_0 \equiv Lie(L) \simeq \frak{o}(2) + \frak o(p-2,q).
$$
Any character of the Lie algebra $\frak l_0$
 (or any complex character of $\frak l$) is determined by
 its restriction to $\hc_0$.
So, we shall write $\Bbb C_\nu$ for 
 the character of the Lie algebra $\frak l_0$
 whose restriction to $\hc$ is $\nu \in (\hc)^*$.
With this notation,
  the character of $L$ acting on 
 $\wedge^{\dim \frak u}\frak u$
 is written as $\Bbb C_{2 \rho(\frak u)}$
 where 
$$
\rho(\frak u)  := (\frac{p+q}2 - 1) f_1.
\tag \num.1
$$

The homogeneous manifold $G/L$ carries a $G$-invariant complex structure
  with canonical bundle $\wedge^{\operatorname{top}} T^* G/L
   \simeq G \times_L \Bbb C_{2 \rho(\frak u)}$.
 As an algebraic analogue of a Dolbeault cohomology 
 of a $G$-equivariant holomorphic vector bundle over a complex manifold $G/L$, 
 Zuckerman introduced the cohomological parabolic induction 
$\Cal R^j_\frak q\equiv\left(\Cal R^\frak g_\frak q\right)^j\,\,(j\in\Bbb N)$, 
which is a covariant functor from the
 category of metaplectic \metalk-modules to 
 that of $(\frak g,K)$-modules.
Here,
 $\widetilde L$ is a metaplectic covering of $L$
 defined by the character of $L$ acting on
 $\wedge^{\dim \frak u}\frak u \simeq \Bbb C_{2 \rho(\frak u)}$.
In this paper,
 we follow the normalization in {\xvr}, Definition 6.20
 which is different from the one in \xvg\ by a \lq$\rho$-shift\rq.

The character $\Bbb C_{\lambda f_1}$ of $\frak l_0$
 lifts to a metaplectic \metalk-module
if and only if $\lambda \in {\Bbb Z} + \frac{p+q}2$.
In particular, we can define \gk-modules
$\zdfc_q^{j}{\lambda f_1}$
 for $\lambda \in A_0(p,q)$.
The $\Cal Z(\frak g)$-infinitesimal character of 
 $\zdfc_q^{j}{\lambda f_1}$
 is given by
$$
     (\lambda,\frac{p+q}2-2,\frac{p+q}2-3,\dots,\frac{p+q}2-[\frac{p+q}2])
    \in (\hc)^*
$$
 in the Harish-Chandra parametrization
 if it is non-zero.
In the sense of Vogan \xvi, 
 we have
$$
\alignat2
&\Bbb C_{\lambda f_1} \text{ is in the good range }
&& \ \Leftrightarrow \
\lambda > \frac{p+q}2-2,
\\
&\Bbb C_{\lambda f_1} \text{ is in the weakly fair range }
&& \ \Leftrightarrow \
\lambda \ge 0.
\endalignat
$$
We note that
$\zdfc_q^{j}{\lambda f_1} = 0$
 if $j \neq p-2$ and if $\lambda \in A_0(p,q)$.
This follows from a general result in \xvu\ for $\lambda \ge 0$;
and \xkupq\ for  $\lambda = -\frac12$.

\def \sc{4}
\sec{}
For $b \in \Bbb Z$,
 we define an algebraic direct sum of 
 $K = O(p) \times O(q)$-modules by 
$$
  \kK{b} \equiv
  \kR{p}{q}{b} :=
     \bigoplus_{\Sb m, n \in \Bbb N
           \\
               m -n \ge b
           \\
               m -n \equiv b \mod 2
            \endSb}
      \spr{m}{p} \boxtimes \spr{n}{q}.
\tag \num.1
$$
For $\lambda \in A_0(p,q)$,
 we put
$$
\alignat1
    b \equiv b(\lambda, p, q) &:= \lambda -\frac p 2 + \frac q 2 +1 \in \Bbb Z,
\tag \num.2
\\
    \epsilon \equiv \epsilon(\lambda, p, q)&:= (-1)^b.
\tag \num.3 
\endalignat
$$

We define the line bundle $\Cal L_n$ over $G/L$ by the character
$nf_1$ of $L$ (see \S 5.3).

Here is a summary for 
 different realizations of the representation $\pip{p}{q}{\lambda}$:
\proclaim{Fact~\num}
Let $p, q \in \Bbb N$ ($p > 1$).
\newline
{\rm 1)}\enspace\
For any $\lambda \in A_0(p,q)$,
 each of the following 5 conditions defines uniquely a \gk-module,
 which are mutually isomorphic.
We shall denote it by $(\pip{p}{q}{\lambda})_K$.
The \gk-module $(\pip{p}{q}{\lambda})_K$ is non-zero and irreducible.
\newline
{\rm i)}\enspace
 A subrepresentation of the degenerate principal representation
 $\princeK{\lambda}{\epsilon}$ (see \S 3.7)
 with $K$-type $\kK{b}$.
\newline
{\rm i)}$'$\enspace
 A quotient of $\princeK{-\lambda}{\epsilon}$
 with $K$-type $\kK{b}$.
\newline
{\rm ii)}\enspace
 A subrepresentation of $C_\lambda^\infty(X(p,q))_K$
 with $K$-type $\kK{b}$.
\newline
{\rm iii)}\enspace
The underlying \gk-module of the Dolbeault cohomology group
\newline
$
     H_{\bar \partial}^{p-2} 
      (G/L, \Cal L_{(\lambda + \frac{p+q-2}{2})})_K.
$
\newline
{\rm iii)}$'$\enspace
The Zuckerman-Vogan derived functor module $\zdfc_{q}^{p-2}{\lambda f_1}$.
\newline
{\rm 2)}\enspace
In the realization of {\rm (iii)},
 if $f \in (\pip{p}{q}{\lambda})_K$, 
 then there exists an analytic function
 $a \in C^{\infty}(S^{p-1} \times S^{q-1})$ such that
$$
     f(\omega \cosh t, \eta \sinh t)
  = a(\omega, \eta)
          e^{-(\lambda + \rho)t}
              (1 + t e^{-2t} O(1))
  \quad  \text{ as } t \to \infty
$$
Here, we put $\rho = \frac{p+q-2}2$.
\endproclaim
For details, we refer, for example,
  to \xhowetan\ for (i) and (i)$'$;
 to \xschlap\ for (ii) and also for a relation with (i)
 (under some parity assumption on eigenspaces);
  to \xkdecomp, \S 6 (see also \xkupq) for (iii)$'$ $\Leftrightarrow$ (ii);
 and to \xwong\ for (iii) $\Leftrightarrow$ (iii)$'$.
 The second statement follows from a general theory of
  the boundary value problem with regular singularities;
  or also follows from a classical asymptotic formula
  of hypergeometric functions (see (8.3.1))
 in our specific setting.
\remark{Remark}
1)\enspace
By definition, (i) and (i)$'$ make sense for $p > 1$ and $q > 0$;
 and others for $p > 1$ and $q \ge 0$.
\newline
2)\enspace
Each of the realization (i), (i)$'$, (ii), and (iii)
 also gives a globalization of $\pip{p}{q}{\lambda}$,
 namely, a continuous representation of $G$ on 
 a topological vector space.
Because all of $(\pip{p}{q}{\lambda})_K$ ($\lambda \in A_0(p,q)$)
 are unitarizable
 we may and do take the globalization $\pip{p}{q}{\lambda}$
 to be the unitary representation of $G$.
\newline
3)\enspace
If $\lambda > 0$ and $\lambda \in A_0(p,q)$, 
 then the realization (ii) of $\pip{p}{q}\lambda$ 
 gives a discrete series representation for $X(p,q)$.
Conversely, 
$$
   \set{\pip{p}{q}{\lambda}}{\lambda \in A_0(p,q), \lambda > 0}
$$
 exhausts the set of discrete series representations for $X(p,q)$.
\endremark
If $(p,q)=(1,0)$,
 then $O(p,q) \simeq O(1)$ and it is convenient to 
 define representations of $O(1)$ by
$$
\pip{1}{0}{\lambda} = 
\cases \boldkey 1 \ & (\lambda = - \frac12),
\\
    \sgn \ & (\lambda = \frac12),
\\
      0 \ & \text{ (otherwise).}
\endcases
$$

As we defined $\pip{p}{q}{\lambda}$ in Fact~\num,
 we can also define an irreducible unitary representation,
 denoted by $\pim{p}{q}{\lambda}$,
 for $\lambda \in A_0(q,p)$ 
 such that the underlying \gk-module has the following
 $K$-type
$$
     \bigoplus_{\Sb m, n \in \Bbb N
           \\
               m -n \le -\lambda +\frac{q}2-\frac{p}2-1
           \\
               m -n \equiv -\lambda +\frac{q}2-\frac{p}2-1 \mod 2
            \endSb}
      \spr{m}{p} \boxtimes \spr{n}{q}.
$$
Similarly to $\pip{p}{q}{\lambda}$,
 the representations $\pim{p}{q}{\lambda}$
 are realized in function spaces on another hyperboloid $O(p,q)/O(p,q-1)$.

In order to understand the notation here,
 we remark:
\item{i)}
$\pim{p}{q}{\lambda} \in \widehat{O(p,q)}$
 corresponds to the representation $\pip{q}{p}{\lambda} \in \widehat{O(q,p)}$
 if we identify $O(p,q)$ with $O(q,p)$. 
\item{ii)}
$
     \pip {p} {0}{\lambda} \simeq \spr{k}{p},
$
where
$
      k=\lambda -\frac {p-2} 2
$
and
$
p \ge 1, k \in \Bbb N.
$

\def\sc{5}
\sec{}
The case $\lambda = \pm\frac{1}{2}$ is delicate,
 which happens when $p+q \in 2 \Bbb N +1$.

First, we assume $p+ q \in 2 \Bbb N + 1$.
By using the equivalent realizations of $\pip{p}{q}{\lambda}$ in Fact~\ch.4
 and by the classification of the composition series of
 the most degenerate principal series representation
 $\princeK{\lambda}{\epsilon}$ (see \xhowetan),
 we have non-splitting short exact sequences of \gk-modules:
$$
\alignat4
    0
    &\to &(\pim{p}{q}{-\frac12})_K
    &\to &\princeK{-\frac12}{(-1)^{\frac{p-q+1}2}}
    &\to &(\pip{p}{q}{\frac12})_K
    &\to 0,
\tag \num.1
\\
    0 
    &\to &(\pip{p}{q}{-\frac12})_K
    &\to &\princeK{-\frac12}{(-1)^{\frac{p-q-1}2}}
    &\to &(\pim{p}{q}{\frac12})_K
    &\to 0.
\tag \num.2
\endalignat
$$
Because $\pip{p}{q}{\lambda}$ ($\lambda \in A_0(p,q)$) is self-dual,
 the dual \gk-modules of (\num.1) and (\num.2) give
 the following non-splitting short exact sequences of \gk-modules:
$$
\alignat4
 0  &\to &(\pip{p}{q}{\frac12})_K
    &\to &\princeK{\frac12}{(-1)^{\frac{p-q+1}2}}
    &\to &(\pim{p}{q}{-\frac12})_K
    &\to 0,
\tag \num.3
\\
 0  &\to &(\pim{p}{q}{\frac12})_K
    &\to &\princeK{\frac12}{(-1)^{\frac{p-q-1}2}}
    &\to &(\pip{p}{q}{-\frac12})_K
    &\to 0.
\tag \num.4
\endalignat
$$

Next, we assume $p+q \in 2 \Bbb N$.
Then, $\varpi^{p,q}$
  is realized as a subrepresentation of some degenerate principal series
 (see \xkorsI, Lemma~3.7.2).
More precisely, 
 we have non-splitting short exact sequences of \gk-modules
$$
\alignat1
 &0   \to   \varpi^{p,q}_K
     \to  \princeK{-1}{(-1)^{\frac{p-q}2}}
     \to  \left((\pim{p}{q}{1})_K \oplus (\pip{p}{q}{1})_K\right)
     \to 0,
\tag \num.5
\\
 &0  
     \to  \left((\pim{p}{q}{1})_K \oplus (\pip{p}{q}{1})_K\right)
     \to  \princeK{1}{(-1)^{\frac{p-q}2}}
     \to   \varpi^{p,q}_K
     \to 0,
\tag \num.6
\endalignat
$$
 and an isomorphism of \gk-modules:
$$
  \princeK{0}{(-1)^{\frac{p-q+2}2}}
  \simeq
  (\pim{p}{q}{0})_K \oplus (\pip{p}{q}{0})_K.
\tag \num.7
$$
These results will be used 
 in another realization of the unipotent representation $\varpi^{p,q}$,
 namely,
 as a submodule of the Dolbeault cohomology group
 in a subsequent paper (cf. Part 1, Introduction, Theorem~B~(4)).

\def \ch{6}
\def \sc{1}
\head
 \S \ch. Conformal embedding of the hyperboloid
\endhead
This section prepares the geometric setup 
 which will be used in \S 7 and \S 9 
 for the branching problem of $\varpi^{p,q}|_{G'}$.
Throughout this section,
 we shall use the following notation: 
$$
\alignat2
    & |x|^2:=|x'|^2 + |x''|^2=\sum_{i=1}^{p'} (x_i')^2 
                            +\sum_{j=1}^{p''} (x_j'')^2,
  &\text{ for }
  &x:= (x',x'') \in \Bbb R^{p' + p''}= \Bbb R^p,
\\
   &  |y|^2:=|y'|^2 + |y''|^2=\sum_{i=1}^{q'} (y_i')^2 
                            +\sum_{j=1}^{q''} (y_j'')^2,  
  &\text{ for }
  &y:= (y',y'') \in \Bbb R^{q' + q''}= \Bbb R^q.
\endalignat
$$
\sec{}
We define two open subsets of $\Bbb R^{p+q}$ by
$$
  \align
     \Bbb R_+^{p'+p'',q'+q''}
    &:=\set{(x,y)=\left((x',x''),(y',y'')\right)\in \Bbb R^{p'+p'',q'+q''}}
         {|x'|>|y'|},
\\
     \Bbb R_-^{p'+p'',q'+q''}
    &:=\set{(x,y)=\left((x',x''),(y',y'')\right)\in \Bbb R^{p'+p'',q'+q''}}
         {|x'|<|y'|}.  
  \endalign
$$
Then
  the disjoint union
 $\Bbb R_+^{p'+p'',q'+q''} \cup \Bbb R_-^{p'+p'',q'+q''}$
 is open dense in $\Bbb R^{p+q}$.
Let us consider the intersection of $\Bbb R_\pm^{p'+p'',q'+q''}$ with
 the submanifolds given in \S 3.2:
$$
   M \subset\Xi \subset \Bbb R^{p,q}.
$$
Then, we define two open subsets of
 $M \simeq S^{p-1} \times S^{q-1}$ by
$$
     M_{\pm} := M \cap \Bbb R_{\pm}^{p'+p'',q'+q''}.  
  \tag{\num.1}
$$
Likewise,
 we define two open subsets of the cone $\Xi$ by
$$
     \Xi_{\pm} := \Xi \cap \Bbb R_{\pm}^{p'+p'',q'+q''}.  
  \tag{\num.2}
$$
 
We notice that if $(x,y) = ((x', x''), (y',y'')) \in \Xi$
 then
$$
    |x'| > |y'| 
    \Longleftrightarrow
    |x''| < |y''| 
$$
because $|x'|^2 + |x''|^2 = |y'|^2 + |y''|^2$.
The following statement is immediate from definition:
$$
  \alignat2
     \Xi_+ = \emptyset
    &\ \Leftrightarrow \
     M_+ = \emptyset
 && \ \Leftrightarrow \quad
   p' q''=0.
\tag \num.3
\\
     \Xi_- = \emptyset
    &\ \Leftrightarrow\
     M_- = \emptyset
 && \ \Leftrightarrow \quad
    p'' q'=0.  
\tag \num.4
  \endalignat
$$
\def\sc{2}
\sec{}
We embed the direct product of hyperboloids
$$
  X(p',q') \times X(q'',p'')
  =\set{((x',y'),(y'',x''))}{|x'|^2-|y'|^2=|y''|^2-|x''|^2=1}.
$$
 into $\Xi_+ \ (\subset \Bbb R^{p,q})$ by the map
$$
  X(p',q') \times X(q'',p'') \hookrightarrow \Xi_+, %
  \ 
   \left((x', y'), (y'',x'')\right) \mapsto (x', x'', y', y'').
\tag \num.1
$$
The image is transversal to rays (see \xkorsI, \S 3.3 for definition)
 and the induced pseudo-Riemannian metric
 $g_{X(p',q') \times X(q'',p'')}$ on
  $X(p',q') \times X(q'',p'')$
 has signature $(p'-1, q') + (p'',q''-1) = (p-1,q-1)$.
With the notation in \S 5.2,
 we have
$$
      g_{X(p',q') \times X(q'',p'')} 
     = g_{X(p',q')} \oplus (-g_{X(q'',p'')}).
$$  
We note that
 if $p'' = q' =0$,
 then $X(p',q') \times X(q'',p'')$ is diffeomorphic to
 $S^{p-1} \times S^{q-1}$,
 and
 $g_{X(p',0) \times X(q'',0)}$
 is nothing but the pseudo-Riemannian metric
 $g_{S^{p-1} \times S^{q-1}}$ of signature $(p-1, q-1)$ (see \xkorsI, \S 3.3).

By the same computation as in (3.4.1),
 we have the relationship among the Yamabe operators on hyperboloids
 (see also (5.2.2)) by
$$
     \tilLap{X(p',q') \times X(q'',p'')} 
=\tilLap{X(p',q')} -\tilLap{X(q',p'')}.  
\tag \num.2
$$

We denote by $\Phi_1$ the composition of (\num.1) and 
 the projection $\Phi\: \Xi \to M$ (see \xkorsI, (3.2.4)),
 namely,
$$
 \Phi_1 \:  X(p',q') \times X(q'',p'') \hookrightarrow M,
 \
   \left( (x',y'),(y'',x'')\right)
          \mapsto \left(\frac{(x',x'')}{|x|}, \frac{(y', y'')}{|y|}\right).
\tag \num.3
$$
\proclaim{Lemma \num}
{\rm{1)}}\enspace
The map  $\Phi_1 \:  X(p',q') \times X(q'',p'') \to M$ is 
 a diffeomorphism onto $M_+$.
The inverse map
 $\Phi_1^{-1} \:  \ M_+  \to X(p',q') \times X(q'',p'')$
 is given by the formula:
$$
    \left( (u',u''),(v',v'')\right)
    \mapsto \left(\frac {(u',v')}{\sqrt{|u'|^2-|v'|^2}},
                   \frac {(v'',u'')}{\sqrt{|v''|^2-|u''|^2}}\right).  
\tag \num.4
$$
\item{\rm{2)}}\enspace
 $\Phi_1$ is a conformal map with conformal factor $|x|^{-1} = |y|^{-1}$,
 where $x=(x',x'') \in \Bbb R^{p'+p''}$
 and $y=(y',y'') \in \Bbb R^{q'+q''}$.  
Namely, we have
$$
     \Phi_1^* \left(g_{S^{p-1} \times S^{q-1}}\right)
    = \frac 1 {|x|^2} 
   g_{X(p',q') \times X(q'',p'')}.  
$$
\endproclaim
\demo{Proof}
The first statement is straightforward
 in light of the formula
$$
     |u'|^2 -|v'|^2 = |v''|^2 -|u''|^2 >0 
$$
 for $(u,v) = ((u', u''), (v', v'')) \in M_+ \subset S^{p-1} \times S^{q-1}$.
\newline
The second statement is a special case of Lemma~3.3.
\qed
\enddemo

\def\sc{3}
\sec{}
Now, 
 the conformal diffeomorphism $\Phi_1
 \: X(p',q') \times X(q'',p'') \rarrowsim M_+$
 establishes a bijection of the kernels of the Yamabe operators
 owing to Proposition~2.6:
\proclaim{Lemma \num}
$\widetilde{\Phi_1^*}$ gives a bijection
 from $\Ker \tilLap{M_+}$
 onto $\Ker \tilLap{X(p',q') \times X(q'',p'')}$.  
\endproclaim
Here,
 the twisted pull-backs $\widetilde{\Phi_1^*}$ and
 $\widetilde{(\Phi_1^{-1})^*}$ 
 (see Definition~2.3),
 namely,
$$
\alignat5
     &\hphantom{i}\widetilde{\Phi_1^*} &&\: 
             && \hphantom{MMMM}C^{\infty}(M_+) 
             &&\;\to\; 
             && C^{\infty}(X(p',q') \times X(q'',p'')),
\tag \num.1
\\
     &\widetilde{{(\Phi_1^{-1}})^*} 
     &&\: &&C^{\infty}(X(p',q')\times X(q'',p''))
                       &&\;\to\; 
            &&\hphantom{MMMM}C^{\infty}(M_+),
\tag \num.2
\endalignat
$$
 are given by the formulae
$$
\alignat1
     (\widetilde{\Phi_1^*} F)(x',y',y'',x'')
    &:=(|x'|^2+|x''|^2)^{-\frac {p+q-4}{4}}
      F\left(\frac {(x',x'')}{\sqrt{|x'|^2+|x''|^2}},
             \frac {(y',y'')}{\sqrt{|y'|^2+|y''|^2}}\right),  
\\
     (\widetilde{(\Phi_1^{-1})^*}f)(u',u'',v',v'')
     &:=(|u'|^2 -|v'|^2)^{-\frac {p+q-4}{4}}
        f\left(\frac {(u',v')}{\sqrt{|u'|^2-|v'|^2}},
          \frac {(u'',v'')}{\sqrt{|v''|^2-|u''|^2}}\right),
\endalignat
$$
respectively.
We remark that
 $\widetilde{(\Phi_1^{-1})^*}=(\widetilde{\Phi_1^*})^{-1}$. 

\def\sc{4}
\sec{}
Similarly to \S \ch.2,
 we consider another embedding
$$
 X(q',p') \times X(p'',q'') \hookrightarrow \Xi_-,
 \quad
   \left( (y',x'),(x'',y'')\right) \mapsto (x', x'', y', y'').
\tag \num.1
$$
The composition of (\num.1) and
 the projection $\Phi \: \Xi \to M$ is denoted by
$$
   \Phi_2 \:  X(q',p') \times X(p'',q'') \hookrightarrow M,
   \
   \left( (y',x'),(x'',y'')\right) 
          \mapsto \left(\frac{(x',x'')}{|x|}, \frac{(y', y'')}{|y|}\right).
\tag \num.2
$$
Obviously, results analogous to Lemma~\ch.2 and Lemma~\ch.3
 hold for $\Phi_2$.
For example, here is a lemma parallel to Lemma~\ch.2:
\proclaim{Lemma \num}\enspace
The map
$
     \Phi_2 \: X(q',p') \times X(p'',q'') \to M_-
$
 is a conformal diffeomorphism onto $M_-$.
The inverse map $\Phi_2^{-1} \:   M_- \to X(q',p') \times X(p'',q'')$ 
 is given by
$$
     \left( (u',u''),(v',v'')\right)
     \mapsto \left(\frac {(v',u')}{\sqrt{|v'|^2-|u'|^2}},
                   \frac {(u'',v'')}{\sqrt{|u''|^2-|v''|^2}}\right).  
$$
\endproclaim

\def \ch{7}
\def \sc{1}
\head
 \S \ch. 
  Explicit branching formula (discrete decomposable case)
\endhead

If one of $p', q', p''$ or $q''$ is zero,
 then
 the restriction $\varpi^{p,q}|_{G'}$
 is decomposed discretely into irreducible representations of
 $G'=O(p',q') \times O(p'',q'')$
 as we saw in \S 4.
In this case,
 we can determine the branching laws of $\varpi^{p,q}|_{G'}$ as follows:
\proclaim{Theorem \num}
Let $p+q \in 2 \Bbb N$.
If $q'' \ge 1$ and $q' + q'' = q$,
 then we have an irreducible decomposition
 of the unitary representation $\varpi^{p,q}$
 when restricted to $O(p,q') \times O(q'')$:
$$
\alignat1
 \varpi^{p,q}|_{O(p, q') \times O(q'')}
 &\simeq
 \overset{\infty\hphantom{M}}\to{\Hsum{l=0\hphantom{M}}}
  \pip{p}{q'}{l + \frac{q''}2 -1} \boxtimes \pim{0}{q''}{l + \frac{q''}2 -1}
\\
 &\simeq
 \overset{\infty\hphantom{M}}\to{\Hsum{l=0\hphantom{M}}}
  \pip{p}{q'}{l + \frac{q''}2 -1} \boxtimes \spr{l}{q''}.
\tag \num.1
\endalignat
$$
\endproclaim
We shall prove Theorem~\num\ in \S \ch.5
 after we prepare an algebraic lemma in \S \ch.3 and 
 a geometric lemma in \S \ch.4.

\remark{Remark}
The formula in Theorem~\ch.1 is nothing but a $K$-type formula
 (see Theorem~3.6.1) when $q'=0$.
\endremark

\def\sc{2}
\sec{}
The branching law (\ch.1.1) is an infinite direct sum for $q'' > 1$.
This subsection treats the case $q'' = 1$,
 which is particularly interesting,
 because the branching formula consists of
 only two irreducible representations
 (we recall $\spr{l}{1} \neq 0$ if and only if $l = 0, 1$).
For simplicity, we shall assume $q \ge 3$ in \S \num.

It follows from Theorem~\ch.1 with $q' =1$ that 
$$
\alignat1
 \varpi^{p,q}|_{O(p, q-1) \times O(1)}
 &\simeq
  \left(\pip{p}{q-1}{-\frac 12} \boxtimes \pim{0}{1}{-\frac12}\right)
  \oplus
  \left(\pip{p}{q-1}{\frac 12} \boxtimes \pim{0}{1}{\frac12}\right)
\\
 &\simeq
  \left(\pip{p}{q-1}{-\frac 12} \boxtimes \boldkey 1\right)
  \oplus
  \left(\pip{p}{q-1}{\frac 12} \boxtimes \sgn\right).
\tag \num.1
\endalignat
$$
This means that $\varpi^{p,q}_K$ can be realized
 in a subspace of $C^\infty_{\frac12}(X(p,q-1))$,
 namely,
 the kernel of the Yamabe operator $\tilLap{X(p,q-1)}$ (see (5.2.3)).

More precisely, according to the direct sum decomposition (see (5.2.4)),
 we have
$$
 \Ker \tilLap{X(p,q-1)}
 \equiv C_{\frac12}^\infty(X(p,q-1)) =
   C_{\frac12}^\infty(X(p,q-1))_+ + 
   C_{\frac12}^\infty(X(p,q-1))_-.
\tag \num.2
$$
We recall that the central element $-I_{p+q} \in G$
 acts on $\varpi^{p,q}$ with scalar $\delta$,
 where 
$$
         \delta := (-1)^{\frac{p-q}2}.
\tag \num.3
$$

In view of the composition series of eigenspaces
 on the hyperboloid
 (see \xschlap\ for the case $\delta = +$; similar for $\delta = -$),
 we have non-splitting exact sequences of Harish-Chandra modules of
 $O(p,q-1)$:
$$
\alignat4
 0  &\to &(\pip{p}{q-1}{-\frac12})_K
    &\to & \left(C_{\frac12}^\infty(X(p,q-1))_\delta\right)_K
    &\to &(\pim{p}{q-1}{\frac12})_K
    &\to 0,
\tag \num.4
\\
 0  &\to &(\pip{p}{q-1}{\frac12})_K
    &\to & \left(C_{\frac12}^\infty(X(p,q-1))_{-\delta}\right)_K  
    &\to &(\pim{p}{q-1}{-\frac12})_K
    &\to 0,
\tag \num.5
\endalignat
$$

Here is a realization of $\varpi^{p,q}$
 in a subspace of 
 the kernel of the Yamabe operator on the hyperboloid $X(p,q-1)$,
 on which $O(p,q)$ acts as meromorphic conformal transformations.
\proclaim{Corollary~\num.1}
Let $W_\pm$ be the unique non-trivial subrepresentation of $O(p,q-1)$
 in $(\Ker \tilLap{X(p,q-1)})_\pm$
 $\equiv$  $C_{\frac12}^\infty(X(p,q-1))_\pm$.
Each of the underlying \gk-module is infinitesimally unitarizable,
 and we denote the resulting unitary representation by $\overline W_\pm$. 
Then,
 the irreducible unitary 
 representation $\varpi^{p,q}$ of $O(p,q)$ is realized on
 the direct sum $\overline{W_+} + \overline{W_-}$.
\endproclaim
We note that $\overline{W_{-\delta}} \subset L^2(X(p,q-1))$
 and $\overline{W_{\delta}} \not\subset L^2(X(p,q-1))$
 where $\delta = (-1)^{\frac{p-q}2}$.

It is interesting to note that
 the Laplacian $\Delta_{X(p,q-1)}$ acts on
 a discrete series $\pip{p}{q-1}{\lambda}$
 for the hyperboloid $X(p,q-1)$
 as a scalar $-\lambda^2 + \frac{1}{4}(p+q-3)^2$,
 which attains the maximum when $\lambda = \frac{1}{2}$
 if $p + (q-1) \in 2 \Bbb N + 1$.

Taking the direct sum of two exact sequences (\num.4) and (\num.5),
 we have the following:
\proclaim{Theorem~\num.2}
There is a non-split exact sequence of Harish-Chandra modules for $O(p-1,q)$.
$$
 0  \to (\varpi^{p,q})_K
    \to (\Ker \tilLap{X(p,q-1)})_K
    \to (\varpi^{p+1, q-1})_K
    \to 0.
\tag \num.6
$$
\endproclaim
It is a  mysterious phenomenon in (\num.6)
 that $\varpi^{p,q}$ extends to a representation of $O(p,q)$
 and $(\varpi^{p+1, q-1})_K$ to that of $O(p+1, q-1)$.
So, different real forms of $O(p+q, \Bbb C)$ act on subquotients
 of
 the kernel of the Yamabe operator on the hyperboloid
 $X(p,q-1) = O(p,q-1)/O(p-1,q-1)$ !
 
Here, we remark that
$\Ker \tilLap{X(p,q-1)} \cap L^2(X(p,q-1)) \neq \{0\}$
 if and only if $p+q \in 2 \Bbb Z$,
  by the classification of discrete series for the hyperboloid $X(p,q-1)$
   for $p>1$.

\def\sc{3}
\sec{}
By Theorem~4.2,
 the restriction $\varpi^{p,q}|_{K'}$ is \adm{K'},
 where $K' \simeq O(p) \times O(q') \times O(q'')$.
Let us first find the $K'$-structure of $\varpi^{p,q}$.
We recall a classical branching law
 with respect to $(O(q), O(q') \times O(q''))$
 ($q = q'+q'', q' \ge 1, q'' \ge 1$):
$$
       \spr{n}{q} 
    \simeq 
       \bigoplus \Sb k, l \in \Bbb N \\ k+l \le n \\ k+l \equiv n \endSb
       \spr{k}{q'} \boxtimes \spr{l}{q''}.
\tag \num.1
$$
Then we have isomorphisms as $K'$-modules:
$$
\alignat1
&\bigoplus
\Sb m, n \in \Bbb N \\ m-n = b \\ m - n \equiv b \mod 2 \endSb
\spr{m}{p} \boxtimes \spr{n}{q}|_{O(p) \times O(q') \times O(q'')}
\\
&\simeq
\bigoplus
\Sb m, n \in \Bbb N \\ m-n = b \\ m - n \equiv b \mod 2 \endSb
\bigoplus
\Sb k, l \in \Bbb N \\ k+l \le n \\ k+l \equiv n \mod 2 \endSb
\spr{m}{p} \boxtimes \spr{k}{q'} \boxtimes \spr{l}{q''}
\\
&\simeq
\bigoplus
\Sb l \in \Bbb N \endSb
\bigoplus
\Sb m, k \in \Bbb N \\ m-b \ge k + l \\ m-b \equiv k+l \mod 2 \endSb
\spr{m}{p} \boxtimes \spr{k}{q'} \boxtimes \spr{l}{q''}
\\
&\simeq
\bigoplus \Sb l \in \Bbb N \endSb
\kR{p}{q'}{b+l} \boxtimes \spr{l}{q''}
\endalignat
$$
In view of Theorem~3.6.1, we have proved:
\proclaim{Lemma~\num}
We have an isomorphism of $K'$-modules:
$$
  \varpi^{p,q}_K \simeq
 \bigoplus \Sb l \in \Bbb N \endSb
 \kR{p}{q'}{\frac{q-p}2+l} \boxtimes \spr{l}{q''}
\tag \num.2
$$
\endproclaim

\def \sc{4}
\sec{}
By Theorem~4.2 and \xkdecoass, Lemma~1.3,
 the underlying \gk-module $\varpi^{p,q}_K$
 is decomposed into an algebraic direct sum of 
 irreducible $(\frak g', K')$-modules:
$$
    \varpi^{p,q}_{K} 
    \simeq \bigoplus_{\tau} m_\tau \tau 
    \simeq \bigoplus_{\tau_1, \tau_2} 
     m_{\tau_1, \tau_2} \tau_1 \boxtimes \tau_2,
\tag \num.1
$$
 where $m_\tau \in \Bbb N$ and
 $\tau$
 runs over irreducible $(\frak g', K')$-modules
 or equivalently,
 $\tau_1$ runs over irreducible $(\frak g'_1, K'_1)$-modules
 (with obvious notation for $G_1' := O(p,q')$)
 and $\tau_2$ runs over irreducible $O(q'')$-modules.
It follows from Lemma~\ch.3 that
 for each $l$
 there exists a $(\frak g'_1, K'_1)$-module $W_l$
 which is a direct sum of irreducible $(\frak g'_1, K'_1)$-module
 such that $W_l$ is isomorphic to $ \kR{p}{q'}{\frac{q-p}2+l}$
 as $K'_1$-modules.
Let us prove that $W_l$ is in fact irreducible
 as a $(\frak g'_1, K'_1)$-module.

\proclaim{Lemma \num}
$W_l$ is realized in a subspace of $C^\infty_\lambda(X(p,q'))$
 with $\lambda = l + \frac{q''}2 -1$.
\endproclaim
\demo{Proof}
In our conformal construction of $\varpi^{p,q}$ in \S 3,
 we recall that each $K$-finite vector
 of $\varpi^{p,q}$ is an analytic function
 satisfying the Yamabe equation on $M \simeq S^{p-1} \times S^{q-1}$.
By using the conformal diffeomorphism $\Phi_1 \: X(p, q') \times S^{q''-1}
 \to M_+$, an open dense subset of $M$ (see (6.1.4) with $p''=0$),
 we can realize $W_l \times \spr{l}{q''}$
 in the space of smooth functions on $X(p,q') \times S^{q''-1}$
 satisfying the Yamabe equation by the following diagram:
$$
\alignat7
  &\Cal A(M)
  &&\underset{\roman{dense}}\to\subset
  &&C^\infty(M)
  &&\underset{\roman{restriction}}\to\hookrightarrow
  &&C^\infty(M_+)
  &&\underset{\Phi_1^*}\to\rarrowsim 
  &&C^\infty(X(p,q') \times S^{q''-1})  
\\
  & \qquad \cup &&&& \qquad\cup &&&& \qquad\cup &&&& \qquad \cup
\\
   W_l \boxtimes \spr{l}{q''} 
  &\subset \varpi^{p,q}_K 
  &&\underset{\roman{dense}}\to\subset
  && \Ker \tilLap{M}
  &&\underset{\roman{restriction}}\to\hookrightarrow
  &&\Ker \tilLap{M_+}
  &&\underset{\Phi_1^*}\to\rarrowsim 
  &&\Ker \tilLap{X(p,q') \times S^{q''-1}}.  
\endalignat
$$
Because
 $\tilLap{X(p,q') \times S^{q''-1}}
 = \tilLap{X(p',q')} - \tilLap{S^{q''-1}}$
 (see (6.2.2))
 acts on $W_l \boxtimes \spr{l}{q''}$ as $0$,
 and because
 $\tilLap{S^{q''-1}}$ acts on $\spr{l}{q''}$ as a scalar
 $\frac{1}{4} - (l + \frac{q''}2 -1)^2$ (see (3.5.1)),
 we conclude that  $\tilLap{X(p',q')}$ acts on $W_l$
 as the same scalar.  
Hence, Lemma is proved.
\qed
\enddemo
\def\sc{5}
\sec{}
Let us complete the proof of Theorem~\ch.1.
It follows from 
 Lemma~\ch.4 together with the $K_1'$-structure of $W_l$ in \S \ch.3
 that $W_l$ is irreducible and isomorphic to
 $(\pip{p}{q'}{l + \frac{q''-2}2})_{K_1'}$
 as $(\frak g_1', K_1')$-module
 (see Fact~5.4).
Therefore, we have an isomorphism of $(\frak g', K')$-modules
$$
 \varpi^{p,q}_K 
 \simeq
 \bigoplus_{l=0}^\infty
  (\pip{p}{q'}{l + \frac{q''}2 -1})_{K_1'} \boxtimes \spr{l}{q''}
 \quad
 \text{(algebraic direct sum)}.
$$
Taking the closure in the Hilbert space,
 we have (\ch.1.1).
Hence Theorem~\ch.1 is proved.
\qed

\redefine\sc{6}
\sec{}
So far,
 we have not used the irreducibility of $\varpi^{p,q}$ in the branching law.
Although the irreducibility of $\varpi^{p,q}$ $(p,q) \neq (2,2)$
 is known \xbz,
 we can give a new and simple proof for it,
 as an application of the branching formulae in Theorem~7.1.

\proclaim{Theorem~\num}
Let $p, q \ge 2$, $p+q \in 2 \Bbb Z$ and $(p,q) \neq (2,2)$.
Then,
$\varpi^{p,q}$ is an irreducible representation of $O(p,q)$.
\endproclaim
\demo{Proof}
Suppose $W \neq \{0\}$ is a closed invariant subspace of the 
 unitary representation $(\varpi^{p,q}, \overline\Vpq)$.
We want to prove $W = \overline\Vpq$.

Without loss of generality ($p$ and $q$ play a symmetric role),
 we may assume $q \ge 3$,
 and fix $q' \ge 1, q'' \ge 2$ such that $q' + q'' = q$.
We write $\overline\Vpq = 
 \overset{\infty\hphantom{M}}\to{\Hsum{l=0\hphantom{M}}} V_l$
 according to the irreducible decomposition (7.1.1)
 of $G' = G_1' \times G_2' := O(p,q') \times O(q'')$.

Because $W$ is non-zero,
 $W$ contains a $K$-type of the form $\spr{a}{p}\boxtimes\spr{b}{q}$,
 which we fix once for all.
In view of the branching formula (7.3.1), 
 $\spr{b}{q}$ contains a non-zero $O(q'')$-fixed vector.
That is, 
$$
 V_0 := \overline{\Vpq}^{G_2'} \supset
 W^{G_2'} \supset {\spr{a}{p}\boxtimes\spr{b}{q}}^{G_2'} \neq \{0\}
$$
Because $G_1' = O(p,q')$ acts on $V_0$ by $\pip{p}{q'}{\frac{q''}2 -1}$
 (see Theorem~\ch.1) which is irreducible,
 and because
$
 V_0 \supset W^{G_2'} 
$
 is stable under the action of $G_1'$,
 we conclude that $V_0 = W^{G_2'}$.
Hence, we have proved
$$
   W \supset V_0.
\tag \num.1
$$

If $W \neq \overline\Vpq$,
 the orthogonal complement $W^\perp \supset V_0$ by the same argument.
This would contradict to $W \cap W^\perp = \{0\}$.
Therefore,
 $W$ must coincides with $\overline\Vpq$.
Hence,  $\varpi^{p,q}$ is irreducible.
\qed
\enddemo

\def \ch{8}
\def \sc{1}
\head
 \S \ch. Inner product on $\varpi^{p,q}$ and the Parseval-Plancherel formula
\endhead

In this section,
 we prove the Parseval-Plancherel type formula
 for our discrete decomposable branching law
 given in Theorem 7.1.  
Our main result in this section is Theorem \ch.5.

\def \sc{1}
\sec{}
The subsections \S \ch.1 $\sim$ \S \ch.3
 review some results and notation from \xkhcrrest\ without proof.  
We set
$$
  \align
  A_+(p,q)&:=A_0(p,q) \cap \set {\lambda \in \Bbb R}{\lambda >0},
\\
\Lambda_\apm(\lambda) 
  &:=\set{(\lambda', \lambda'') \in A_+(p',q') \times A_+(q'',p'')}{
          \lambda' - \lambda'' - \lambda - 1 \in 2 \Bbb N},
\\
\Lambda_\app(\lambda) 
  &:=\set{(\lambda', \lambda'') \in A_+(p',q') \times A_+(p'',q'')}{
         \lambda - \lambda' - \lambda'' - 1 \in 2 \Bbb N}.  
  \endalign
$$
We recall the Jacobi function
$$
     \varphi_{i \lambda}^{(\lambda'', \lambda')} (t)
   := {}_2F_1\left(
       \frac{\lambda'+\lambda''+1-\lambda}2,
       \frac{\lambda'+\lambda''+1+\lambda}2;
       \lambda''+1;-\sinh^2 t
     \right),
$$
which satisfies
$$
     \varphi_{i \lambda}^{(\lambda'', \lambda')} (0)=1,   
$$
$$
   \left( 
t(1-t)\frac{d^2}{d t^2}
+ (c-(a+b+1)t) \frac{d}{d t}
 - a b \right) 
 \varphi_{i \lambda}(t) = 0
$$
with
$$
  a:= \frac{\lambda'+\lambda''+1-\lambda}2,\quad
  b:= \frac{\lambda'+\lambda''+1+\lambda}2,\quad
  c:=\lambda''+1.
$$
We set meromorphic functions of three variables
 $\lambda$, $\lambda'$ and $\lambda''$ by 
$$
\alignat1
     \V{\apm}{\lambda'}{\lambda''}{\lambda}
     &:=
     \frac{
        (\Gamma(\lambda''+1))^2 
         \ 
         \Gamma(\frac{\lambda' - \lambda''+ \lambda+1}{2})
         \ 
         \Gamma(\frac{\lambda' - \lambda'' - \lambda+1}{2})
          }{
         2 \lambda
         \
         \Gamma(\frac{\lambda' + \lambda''+ \lambda+1}{2})
         \ 
         \Gamma(\frac{\lambda' + \lambda'' - \lambda+1}{2})
         }, 
\\
     \V{\app}{\lambda'}{\lambda''}{\lambda}
     &:=
     \frac{
        (\Gamma(\lambda''+1))^2 
         \ 
         \Gamma(\frac{-\lambda' - \lambda''+ \lambda+1}{2})
         \ 
         \Gamma(\frac{\lambda' - \lambda'' + \lambda+1}{2})
          }{
         2 \lambda
         \
         \Gamma(\frac{-\lambda' + \lambda''+ \lambda+1}{2})
         \ 
         \Gamma(\frac{\lambda' + \lambda'' + \lambda+1}{2})
         }.  
\endalignat
$$
Then we have
$$
  \alignat 2
   &\int_0^\infty |\varphi_{i \lambda}^{(\lambda'', \lambda')}(t)|^2
       (\cosh t)^{2 \lambda'+1}
       (\sinh t)^{2 \lambda''+1}
   \ d t
   = \V\apm{\lambda'}{\lambda''}{\lambda}
  \quad&&\text{ for }
   (\lambda',\lambda'')\in \Lambda_\apm(\lambda),
\\
  &
     \int_0^{\frac {\pi} 2} |\varphi_{i \lambda}^{(\lambda'', \lambda')}(i \theta)|^2
       (\cos \theta)^{2 \lambda'+1}
       (\sin \theta)^{2 \lambda''+1}
   \ d \theta
   = \V\app{\lambda'}{\lambda''}{\lambda}
  &&\text{ for }
   (\lambda',\lambda'')\in \Lambda_\app(\lambda).
  \endalignat
$$

\def\sc{2}
\sec{}
We put
$$
    M\equiv M(\lambda,\lambda',\lambda'')
     :=(-1)^{\frac {\lambda'-\lambda''-\lambda-1}{2}}
       \frac{\Gamma(\frac{\lambda'+\lambda''-\lambda+1}{2})\Gamma(\lambda+1)}
            {\Gamma(\frac{\lambda'-\lambda''+\lambda+1}{2})
            \Gamma(\lambda''+1)}.
\tag{\num.1}  
$$
\proclaim{Lemma \num\ {\rm (triangular relation of the Jacobi function)}}
If $\cot \theta = \cosh t$
($0 < \theta < \frac {\pi}{2}, 0<t$),
then
$$
     \varphi_{i \lambda'}^{(\lambda,\lambda'')}(i \theta)
    =M(\cosh t)^{\lambda+\lambda'+\lambda''+1}
     \varphi_{i \lambda}^{(\lambda'',\lambda')}(t).  
$$
\endproclaim
Furthermore,
 we have
$$
  \lambda' \V{\app}{\lambda''}{\lambda}{\lambda'}
  =M^2 \lambda \V{\apm}{\lambda'}{\lambda''}{\lambda}.  
  \tag{\num.2}
$$

\def\sc{3}
\sec{}
$
  \text{Unitarizaition of $\pip{p}{q}{\lambda}$
 $(\lambda \in A_0(p,q))$.}
$
\par
In this subsection,
 we give an explicit unitary inner product
 on $\pip p q {\lambda}$
 for $\lambda \in A_0(p,q)$.  
In view of Fact 5.4,
 the \gk-module $(\pip p q {\lambda})_K$
 is realized in $C_{\lambda}^{\infty}(X(p,q))_K$
 with $K$-types $\kK b$
 where $b=\lambda - \frac p 2 + \frac q 2 +1$.  
Suppose $f \in C_{\lambda}^{\infty}(X(p,q))_K$ belongs
 to the $K$-type $\spr m p \boxtimes \spr n q \in \kK b$.  
Then $f$ is of the form
$$
 f(\omega \cosh t, \eta \sinh t)
  = h_m(\omega) h_n(\eta) (\cosh t)^m (\sinh t)^n 
   \varphi_{i \lambda}^{(n + \frac{q}2-1, m + \frac{p}2-1)}(t),
  \tag{\num.1}
$$
 where $h_m \in \spr{m}{p}$, $h_n \in \spr{n}{q}$,
 $\omega \in S^{p-1}, \eta \in S^{q-1}$, and $t > 0$.  
We put
$$
  \| f\|_{\pip{p}{q}{\lambda}}^2 
  :=
  \| h_m\|_{L^2(S^{p-1})}^2 
  \
  \| h_n\|_{L^2(S^{q-1})}^2 
  \
   \lambda  \V{\apm}{m + \frac{p-2}{2}}{n + \frac{q-2}2}{\lambda}.  
  \tag{\num.2}
$$
This norm defines an inner product
 on $(\pip p q {\lambda})_K$
 ($\lambda \in A_0(p,q)$)
 and one defines an irreducible unitary representation
 of $G=O(p,q)$
 on its Hilbert completion.  
There is an obvious inner product for $\lambda >0$,
 because $(\pip p q {\lambda})_K \subset L^2(X(p,q))$.  
The relation between our norm
 $\|\;\|_{\pip p q {\lambda}}$ 
 and $L^2$-norm $\|\;\|_{L^2(X(p,q))}$
 is given by 
$$
    \| f\|_{\pip{p}{q}{\lambda}}^2  = 
  \lambda \| f\|_{L^2(X(p,q))}^2,
  \quad
  \text{ for any } f \in (\pip{p}{q}{\lambda})_K,   
  \tag{\num.3}
$$
if $\lambda >0$.
For  $\lambda \in A_0(p,q)$
 such that $\lambda \le 0$
 (this can happen if $\lambda =0, -\frac 1 2$),
 $\pip p q {\lambda}$ is also unitarizable
 by this inner product $(\, , \,)_{\pip p q{\lambda}}$
 given in (\num.2),
 which is still positive definite.

\remark{Remark~\num\ (Unitarity)}
As we explained,
 all of $(\pip p q {\lambda})_K$ are unitarizable
 for $\lambda \in A_0(p,q)$.  
We explain four different approaches:  
\item{1)}
If $\lambda > 0$, then $(\pip{p}{q}\lambda)_K$ is unitarizable
 because of the realization in $L^2(X(p,q))$.
\item{2)}
If $\lambda \ge 0$, then $(\pip{p}{q}\lambda)_K$ is unitarizable
 because of the realization of Zuckerman-Vogan's derived
 functor module $\zdfc_q^{p-2}{\lambda}$
 with the parameter $\lambda$ in the weakly fair range \xvu.

We note that the case 
 $\lambda = -\frac{1}2$ is not treated in the above two methods.
However,
\item{3)}
Use the classification of unitarizable subquotients
 of $\princeK{\lambda}{\epsilon}$ given in \xhowetan, \S 3.

At last,
 here is a new proof of the unitarizability
 of $(\pip{p}{q}{\lambda})_K$
 for all $\lambda \in A_0(p,q)$.  
The idea is to use our branching formula Theorem 7.1,
 for which the proof does not use the unitarizability
 of $(\pip p q {\lambda})_K$:
\item{4)}
All of $(\pip{p}{q}{\lambda})_K$ are unitarizable
 because they appear as discrete spectra
 in the branching law of a unitary representation $\varpi^{p,q+c}$
 of a larger group $O(p, q+c)$
 to $O(p,q) \times O(c)$
 for some $c > 0$ (see Theorem~7.1).
For this purpose, 
 $c \le 3$ will do.
\endremark

\def\sc{4}
\sec{}
We notice that the map $\Phi_1 \: X(p,q') \times S^{q''-1} \to M_+$
 (6.2.3) 
 (in the case $p'' = 0$) is given by
$$
 ((\omega \cosh t, \eta' \sinh t), \eta'') \mapsto
  (\omega, (\eta'\cos \theta, \eta'' \sin \theta)),
$$
 where $\theta$ and $t$ satisfy $\cot \theta = \cosh t$.
Suppose $f \in C^{\infty}(X(p,q') \times S^{q''-1})$
 belongs to 
$$
     \pip {p}{q'}{l + \frac {q''}{2}-1} 
    \boxtimes
     \spr l {q''}
$$
 as an $O(p')\times O(q'')$-module, 
 and furthermore to $\spr m p \boxtimes \spr n {q'}$
 as an $O(p) \times O(q')$-module
 in the first factor.  
Then $f$ is of the form:
$$
 f((\omega \cosh t, \eta' \sinh t), \eta'')
  = h_m(\omega) \ h_k(\eta') \ h_l(\eta'')
  (\cosh t)^m (\sinh t)^k
   \varphi_{i (l+\frac{q''}2-1)}^{(k + \frac{q'}2-1, m+\frac{p}{2}-1)}(t)
$$
 where $h_m \in \spr{m}{p}$, $h_k \in \spr{k}{q'}$,
 $h_l \in \spr{l}{q''}$,
 $\omega \in S^{p-1}, \eta' \in S^{q'-1}, \eta'' \in S^{q''-1}$,
 and $t > 0$.

\proclaim{Lemma \num}
The twisted pull-back
 $\widetilde{{(\Phi_1^{-1}})^*} 
     \: C^{\infty}(X(p,q')\times S^{q''-1})
                       \to C^{\infty}(M_+)
$
 (see (6.3.2) for definition) is given by the formula:
$$
\multline
  (\widetilde{{(\Phi_1^{-1}})^*} f)
  (\omega, (\eta'\cos \theta, \eta'' \sin \theta))
\\
  = M^{-1} h_m(\omega) \ h_k(\eta') \ h_l(\eta'') \
    (\cos \theta)^k \ (\sin \theta)^l \
   \varphi_{i(n+\frac{q}{2}-1)}^{(l+\frac{q''}2-1, k + \frac{q'}2-1)}(i \theta).
\endmultline
\tag \num.1
$$
\endproclaim

\demo{Proof}
In view of the definition,
 we have
$$
  (\widetilde{{(\Phi_1^{-1}})^*} f)
  (\omega, (\eta'\cos \theta, \eta'' \sin \theta))
 = (\cosh t)^{\frac{p+q-4}2} 
 f((\omega \cosh t, \eta' \sinh t), \eta'').
$$
Then (\num.1) follows from Lemma~\ch.2 and from $m+ \frac{p}2 = n+\frac{q}2$.
\qed
\enddemo

\def\sc{5}
\sec{}
Let $G' = O(p, q') \times O(q'')$.
We define an inner product on
 the irreducible constituent 
 $\pip{p}{q'}{\lambda} \boxtimes \pim{0}{q''}{\lambda}$
 in the branching law $\varpi^{p,q}|_{G'}$
 (see Theorem~7.1)
 for each $\lambda = l + \frac{q''}2-1$ ($l \in \Bbb N$) by
$$
 \| f \| :=  \| f_1 \|_{\pip{p}{q'}{\lambda}} \ \| f_2 \|_{L^2(S^{q''-1})}
$$
 for $f = f_1 \ f_2$ with $f_1 \in \pip{p}{q'}{\lambda}$ and
 $f_2 \in \pim{0}{q''}{\lambda} \simeq \spr{l}{q''}$.

Because we use an explicit map to prove the branching law,
 the generalized Parseval-Plancherel formula makes sense.
\proclaim{Theorem~\num}
{\rm 1)}\enspace
If we develop $F \in \Ker{\tilLap{M}}$ as
 $F = \sum_l^\infty F_{l}^{(1)} F_{l}^{(2)}$
 according to the irreducible decomposition
$$
  \widetilde{{(\Phi_1})^*}
 \:
 \varpi^{p,q}|_{O(p, q') \times O(q'')}
 \rarrowsim
 \overset{\infty\hphantom{M}}\to{\Hsum{l=0\hphantom{M}}}
  \pip{p}{q'}{l + \frac{q''}2 -1} \boxtimes \spr{l}{q''}
\tag \num.1
$$
 then we have
$$
  \|F\|_{\varpi^{p,q}}^2 
 = \sum_{l=0}^\infty  
  \| F_{l}^{(1)}\|_{\pip{p}{q'}{l+\frac{q''}2-1}}^2 \| 
  \ 
  F_{l}^{(2)}\|_{L^2(S^{q''-1})}^2.
\tag \num.2
$$
\noindent
{\rm 2)}\enspace
In particular,
 if $q'' \ge 3$,
 then all of $\pip{p}{q'}{l + \frac{q''}2-1}$ are
 discrete series for the hyperboloid $X(p,q')$
 and
$$
  \|F\|_{\varpi^{p,q}}^2 
 = \sum_{l=0}^\infty  
  (l + \frac{q''}2-1)
  \ 
  \| F_{l}^{(1)}\|_{L^2(X(p,q'))}^2 \| 
  \ 
  F_{l}^{(2)}\|_{L^2(S^{q''-1})}^2.
\tag \num.3
$$
\endproclaim
\remark{Remark}
The formula (\num.3) coincides with the Kostant-Binegar-Zierau formula 
 (see \S 3)
 in the special case where $q' = 0$
 (namely, where $G'$ is compact).
\endremark
\demo{Proof}
We write $\lambda := l + \frac{q''}2-1$, 
         $\lambda' := m + \frac{p}2-1=n+\frac q 2 -1$
         and
         $\lambda'' := k + \frac{q'}2-1$.
If $F$ is of the form of the right side of (\ch.4.1), then
$$
  \align
\frac{\| \widetilde{{(\Phi_1})^*} F \|^2_{
      \pip{p}{q'}{\lambda} \boxtimes \spr{l}{q''}}
 }{\|F\|_{\varpi^{p,q}}^2}
=&
 \frac{M^2 \ 
 \|h_m\|^2_{L^2(S^{p-1})}
 \ \|h_k\|^2_{L^2(S^{q'-1})} 
 \ \|h_l\|^2_{L^2(S^{q''-1})}
 \ \lambda \V{\apm}{\lambda'}{\lambda''}{\lambda}
      }{
 \|h_m\|^2_{L^2(S^{p-1})}
 \ \|h_k\|^2_{L^2(S^{q'-1})} 
 \ \|h_l\|^2_{L^2(S^{q''-1})}
 \       \lambda' \V{\app}{\lambda''}{\lambda}{\lambda'}
      }
\\
=&1
  \endalign
$$
because of (\ch.2.1).
Hence the first statement is proved.    
The second statement follows from (\ch.3.2).  
\qed
\enddemo

\def \ch{9}
\def \sc{1}
\head
 \S \ch. Construction of discrete spectra in the branching laws
\endhead

\sec{}
In \S 7, we determined explicitly the branching law $\varpi^{p,q}|_{G'}$
 of the minimal unipotent representation $\varpi^{p,q} \in \widehat{G}$
 where $G' = O(p,q') \times O(q'')$.
The resulting branching law has no continuous spectrum 
 (see Theorem~4.2 and Theorem~7.1).
In this section,
 we treat a more general case where continuous spectrum may appear,
 namely,
 the branching law with respect to the semisimple symmetric pair 
$$
(G, G')
= (O(p,q), O(p',q') \times O(p'', q'')),
$$
 where $p' + p'' =p \ (\ge 2)$, $q' + q'' = q \ (\ge 2)$,
 $p+q \in 2 \Bbb N$,  and $(p,q) \neq (2,2)$.

We shall construct explicitly discrete spectra 
 by using the conformal geometry.

Retain the notation in \S 5.1.
We set 
$$
   \Adisc(p,q) := A_0(p,q) \cap \set{\lambda \in \Bbb R}{\lambda >1}.
\tag \num.1
$$
\proclaim{Theorem \num}
The restriction of the unitary representation $\varpi^{p,q}|_{G'}$ contains
$$
\Hsum{\lambda \in \Adisc(p', q') \cap \Adisc(q'',p'')}
  \pip{p'}{q'}{\lambda} \boxtimes \pim{p''}{q''}{\lambda}
\oplus
\Hsum{\lambda \in \Adisc(q', p') \cap \Adisc(p'',q'')}
  \pim{p'}{q'}{\lambda} \boxtimes \pip{p''}{q''}{\lambda}
$$
as a discrete summand.
Here,
  $\pip{p'}{q'}{\lambda} \boxtimes \pim{p''}{q''}{\lambda} \in \widehat{G'}$
 is the outer tensor product of 
  $\pip{p'}{q'}{\lambda} \in \widehat{O(p', q')}$
  and $\pim{p''}{q''}{\lambda} \in \widehat{O(q'',p'')}$.
\endproclaim
We have already established the full branching law $\varpi^{p,q}|_{G'}$
 if $p' q' p'' q'' = 0$.
Thus,
 the main part of this section will be devoted to the proof
 of Theorem~\ch.1 when $p' q' p'' q'' \neq 0$.
We shall give some remarks of Theorem~\num\ at the end of this subsection.
\def\sc{2}
\sec{}
Let $K'= O(p') \times O (q') \times O(p'') \times O(q'')$.
We realize the unitary representation $\varpi^{p,q}$
 on the Hilbert space $\overline{V^{p,q}}$.
Here we recall the notation of \S 3.9, briefly as follows:
$$
\alignat1
    &\overline{V^{p,q}} = \Ker(D_p-D_q)
 \quad \underset{\roman{closed}}\to\subset
         \Cal V = \Dom(D_p) \cap \Dom(D_q)  
 \quad (\underset{\roman{dense}}\to\subset L^2(M))
\\
    & M = S^{p-1} \times S^{q-1} 
 \quad \underset{\roman{dense}}\to\supset \ M_+ \cup M_-.
\endalignat
$$ 
Different from Corollary~4.3 in the discretely decomposable case
 ($p' q' p'' q'' = 0$),
 a $K'$-finite vector of a $G'$-irreducible summand
 (i.e\. a discrete spectrum) in $\varpi^{p,q}|_{G'}$ is
 not necessarily a real analytic function on $M$
 if $p' q' p'' q'' \neq 0$ in our conformal construction of $\varpi^{p,q}$.
With this in mind,
 we extend $\widetilde{(\Phi_1^{-1})^*}f \in C^{\infty}(M_+)$
 (see \S 6 for notation)
 to a function, denoted by $T_+ f$,
 on $S^{p-1} \times S^{q-1}$
 as
$$
     T_+ f := \cases
                  \widetilde{(\Phi_1^{-1})^*} f
                  \quad 
                  &\text{ on } M_+,
\\
                  0                                   
                  &\text{ on } M \setminus M_+. 
            \endcases
\tag \num.1
$$
Here is a key lemma:
\proclaim{Lemma \num}
Suppose $\lambda' \in A_0(p',q')$ and $\lambda'' \in A_0(q'',p'')$.
Let $f \in C^{\infty} (X(p',q') \times X(q'', p''))$ be a $K'$-finite function
 which belongs to $\pip{p'}{q'}{\lambda'} \boxtimes \pim{p''}{q''}{\lambda''}$
 (see Fact~5.4~(1)~(ii)).
\newline
{\rm 1)}\enspace
If $\lambda' \ge \frac{1}{2}$ and $\lambda'' \ge \frac{1}{2}$,
 then $T_+ f \in \Cal V$.
\newline
{\rm 2)}\enspace
If $\lambda' > 1$ and $\lambda'' > 1$,
 then $Y (T_+ f) \in L^2(M)$ and $Y Y' (T_+f) \in L^1(M)$
 for any smooth vector field $Y, Y'$ on $M$.
\newline
{\rm 3)}\enspace
If $\lambda' = \lambda'' > 1$
 then $T_+ f \in \overline{V^{p,q}}$.
\endproclaim

Before proving Lemma~\num,
 we first show that Lemma~\num\ implies Theorem~\ch.1.
In fact, 
 Lemma~\num~(3) constructs a non-zero $(\frak g', K')$-homomorphism
$$
T_+ \:(\pip{p'}{q'}{\lambda})_{K_1'} \boxtimes
      (\pim{p''}{q''}{\lambda})_{K_2'}
     \to \varpi^{p,q},
$$
 for $\lambda \in \Adisc(p', q') \cap \Adisc(q'', p'')$.
Then $T_+$  is injective because
$(\pip{p'}{q'}{\lambda})_{K_1'} \boxtimes
      (\pim{p''}{q''}{\lambda})_{K_2'}$ is irreducible.
Since $\varpi^{p,q}$ is a unitary representation of $G$,
 $T_+$ extends to an isometry of unitary representations of $G'$:
$$
 \pip{p'}{q'}{\lambda} \boxtimes \pim{p''}{q''}{\lambda}
     \to \varpi^{p,q},
$$
 by taking the closure
 with respect to the inner product induced from $\varpi^{p,q}$.
This proves Theorem~\ch.1 for the irreducible representation
$\pip{p'}{q'}{\lambda} \boxtimes \pim{p''}{q''}{\lambda}$
 that appears in the first summand.
The second summand is constructed similarly by using the
conformal diffeomorphism (see Lemma~6.4)
$$
\Phi_2 \:    X(q', p') \times X(p'', q'') \rarrowsim M_- \ (\subset M).
$$
Hence, the proof of Theorem~\ch.1 is completed by assuming Lemma~\num.

\def\sc{3}
\sec{}
\remark{Remark~\num}
There exist
 $\lambda' \in A_0(p',q')$ and $\lambda'' \in A_0(q'',p'')$ 
 satisfying $\lambda' = \lambda''>1$ 
 if and only if $p' + q' \equiv p'' + q'' \mod 2$ and $p' \ge 2, q''\ge 2$.
This implies $p+q \in 2\Bbb N$, $p \ge 2$ and $q\ge 2$,
 and $V^{p,q}$ is non-zero (see \S 3).
Of course,  Lemma~\ch.2~(3) also implies $V^{p,q} \neq \{0\}$. 
\endremark

\def \sc{4}
\sec{}
It follows from Lemma~3.8.1 that
 the first statement of Lemma~\ch.2 is proved if
 $Y (T_+) \in L^{2-\epsilon}$ for any $\epsilon >0$
 and for any smooth vector field $Y$ on $M$.
Therefore,
 both of (1) and (2) of Lemma~\ch.2 are proved
 by studying the asymptotic behaviour of $T_+ f$
 near the boundary of $M_+$ in $M$.
This asymptotic estimate is studied below.

Any $K'$-finite vector $f$ is a finite linear combination
 of the form
 $f_1 f_2$,
 where 
$$
\left\{
\aligned
&\text{$f_1 \in C^\infty(X(p', q'))$ \ is an $O(p') \times O(q')$-finite vector
 that belongs to $\pip{p'}{q'}{\lambda'}$},
\\
&\text{$f_2 \in C^\infty(X(q'', p''))$ is an $O(p'') \times O(q'')$-finite vector
 that belongs to $\pim{p''}{q''}{\lambda''}$.}
\endaligned
\right.
$$  
In order to prove Lemma~\ch.2,
 we may and do assume $f$ is of the form $f_1 f_2$.
Then we have
$$
\multline
           (T_+ (f_1 f_2))(u,v)
\\
          =(|u'|^2 - |v'|^2)_+^{-\frac {p+q-4}{4}}
            f_1\left(\frac{(u',v')}{\sqrt{|u'|^2 - |v'|^2}}\right)
            f_2\left(\frac{(u'',v'')}{\sqrt{|u''|^2 - |v''|^2}}\right).
\endmultline
\tag \num.1
$$
Here,
 we have used the following notation:
$$
    r_+^\nu := \cases r^\nu & (r >0), \\ 0 & (r \le 0). \endcases
$$
\def\sc{5}
\sec{}
In order to analyze the asymptotic behaviour of $T_+ (f_1 f_2)$ 
 (see (\ch.4.1)) near the boundary of $M_+$,
 we consider a change of variables on $S^{p-1} \times S^{q-1}$
 by the surjective map
$$
  \alignat3
           (S^1 \times S^{p'-1}\times S^{p''-1})
    &\times (S^1 \times S^{q'-1} \times S^{q''-1}) 
    &&\to
           S^{p-1} &&\times S^{q-1},
\\
            (e^{i \theta},\omega',\omega'')&,(e^{i \varphi}, \eta',\eta'') 
    &&\mapsto 
            (u,&& v)
  \endalignat
$$
 defined by
$$
     (u,v) \equiv (u',u'',v',v'')
          :=(\omega' \cos \theta, \omega'' \sin \theta,
            \eta' \cos \varphi, \eta'' \sin \varphi).  
\tag \num.1
$$
Because
$
     |u'|^2 -|v'|^2
    =|\omega' \cos \theta|^2 - |\eta' \cos \varphi|^2
    =|\cos \theta|^2 - |\cos \varphi|^2,
$
 $M_\pm$ defined in (6.1.1) is rewritten as 
$$
  \align
     M_+&=\set{(\omega' \cos \theta, \omega'' \sin \theta,
               \eta' \cos \varphi, \eta'' \sin \varphi)}
              {|\cos \theta| > |\cos \varphi|}
\\
     M_-&=\set{(\omega' \cos \theta, \omega'' \sin \theta,
               \eta' \cos \varphi, \eta'' \sin \varphi)}
              {|\cos \theta| < |\cos \varphi|}
  \endalign
$$
Here is an elementary computation
 corresponding to the change of variables (\num.1):
\proclaim{Lemma~\num}
{\rm 1)}\enspace
The volume element $d u\ d v$ on $S^{p-1} \times S^{q-1}$
 is given by
$$
 d u \ d v 
 = |\cos \theta|^{p'-1} \ |\sin \theta|^{p''-1} \
 |\cos \varphi|^{q'-1} \ |\sin \varphi|^{q''-1} \
 d \theta \ d \varphi \ d \omega' \ d \omega'' \ d \eta' \ d \eta'',
$$
 where $d \omega'$ is the volume element on $S^{p'-1}$ and so on.
\newline
{\rm 2)}\enspace
Any smooth vector field on $S^{p-1} \times S^{q-1}$
 is a linear combination of
$$
  \frac{1}{\cos \theta} X', 
  \frac{1}{\sin \theta} X'',
  \frac{1}{\cos \varphi} Y', 
  \frac{1}{\sin \varphi} Y'',
  \der{\varphi}, 
  \der{\theta}, 
$$
  whose coefficients are smooth functions of 
 $(\omega',\omega'', \eta',\eta'', \theta, \varphi)$.
Here,
 $X', X'', Y'$ and $Y''$ are smooth vector fields on 
 $S^{p'-1}, S^{p''-1}, S^{q'-1}$ and $S^{q''-1}$,
 respectively.
\endproclaim
Fact~5.4~(2) describes the asymptotic behaviour
 of $K$-finite functions that belong to 
 discrete series representations for a hyperboloid.
Applying it to $f_1$ with respect to the coordinate (\num.1),
 we find $a_1 \in C^{\infty}(S^{p'-1} \times S^{q'-1})$
 and $h_1 \in C^{\infty}(\Bbb R)$ such that
$$
     f_1 \left(\frac {(u',v')}{\sqrt{|u'|^2 -|v'|^2}}\right)
    =a_1(\omega',\eta')
     \left(\frac {\cos^2 \theta + \cos^2 \varphi}
                 {\cos^2 \theta - \cos^2 \varphi}\right)_+^
     {-\frac {2 \lambda' + p' +q' -2}{4}}
     h_1\left(\sqrt{\frac {\cos^2 \theta - \cos^2 \varphi}
                          {\cos^2 \theta + \cos^2 \varphi}}\right)
$$
Likewise,
 there exist $a_2 \in C^{\infty}(S^{q''-1} \times S^{p''-1})$
 and $h_2 \in C^{\infty}(\Bbb R)$ such that
$$
     f_2 \left(\frac {(v'',u'')}{\sqrt{|v''|^2 -|u''|^2}}\right)
    = a_2(\eta'',\omega'')
     \left(\frac {\sin^2 \theta + \sin^2 \varphi}
                 {\cos^2 \theta - \cos^2 \varphi}\right)_+^
     {-\frac {2 \lambda'' + p'' +q'' -2}{4}}
      h_2\left(\sqrt{\frac {\cos^2 \theta - \cos^2 \varphi}
                          {\sin^2 \theta + \sin^2 \varphi}}\right)
$$
We treat the boundary $\partial M_+$ of $M_+ \subset M$
 locally in the following three cases:
\newline\noindent
{Case 1)}\enspace
 $\cos^2 \theta-\cos^2 \varphi=0$, $(\cos \theta, \cos \varphi) \ne (0,0)$,
 $(\sin\theta,\sin\varphi) \ne (0, 0)$.  
\newline\noindent
{Case 2)}\enspace
 $\cos \theta=\cos \varphi=0$.  
\newline\noindent
{Case 3)}\enspace
 $\sin \theta = \sin \varphi=0$.

We note that Case~(2) or (3) happens only when $p' q' p'' q'' \neq 0$.
\def \sc{6}
\sec{}
(Case 1):\enspace
In this subsection,
 we consider a generic part of the boundary $\partial M_+$
 corresponding to Case 1.
In a local coordinate 
 $(\omega',\omega'', \eta',\eta'', \theta, \varphi)
 \in S^{p'-1}\times S^{p''-1} \times S^{q'-1} \times S^{q''-1}
 \times \Bbb R^2$,
 $T_+ (f_1 f_2))(u,v)$ is written as
$$
 A \ (\cos ^2\theta - \cos ^2 \varphi)_+ 
     ^{-\frac {p+q-4} 4 + \frac {2 \lambda' + p' + q'-2} 4 
       + \frac {2 \lambda'' + p'' + q''-2} 4}
   =
     A \ (\cos ^2\theta - \cos ^2 \varphi)_+ 
     ^{\frac {\lambda' +\lambda''}{2}}.  
$$
Here, $A$ is a smooth function of variables
 $\omega',\omega'', \eta',\eta'', 
 (\cos^2\theta - \cos^2 \varphi)_+^{\frac12}$.
Therefore,
 by using Lemma~\ch.5,
 we have:
$$
\alignat3
&\lambda' + \lambda'' > -1
 &\quad \Rightarrow &\quad T_+ f \in L^2_{\roman loc},&&
\\
&\lambda' + \lambda'' \ge 1
 &\quad \Rightarrow &\quad 
 Y(T_+ f) 
    \in L^{2-\epsilon}_{\roman loc} 
&& \qquad 
\text{for any } \epsilon > 0, Y \in \frak X(M),
\\
&\lambda' + \lambda'' > 2 
&\quad \Rightarrow &\quad
 Y_1 Y_2(T_+ f)    \in L^{1}_{\roman loc}
&& \qquad 
\text{for any } Y_1, Y_2 \in \frak X(M),
\endalignat
$$
 in a neighbourhood of the boundary point of $M_+$ for Case (1).
\def \sc{7}
\sec{}
(Case 2):\enspace
In this subsection
 we consider a neighbourhood of 
 $(\omega',\omega'', \eta',\eta'', \theta, \varphi)$
 satisfying the condition of Case 2.
We take a polar coordinate for
 $(\cos \theta, \cos \varphi) \fallingdotseq (0,0)$ as
$$
  \cases
        \cos \theta &= r \cos \psi,
\\
        \cos \varphi&= r \sin \psi.
  \endcases
$$
The composition to (\ch.5.1) yields a new
 coordinate on $S^{p-1} \times S^{q-1}$ given by
$$
  (u,v)= (u', u'', v', v'')
       = (\omega' r \cos \psi, \omega'' \sqrt{1 - r^2 \cos^2 \psi},
           \eta' r \sin \psi, \eta''    \sqrt{1 - r^2 \sin^2 \psi}),
$$
where
 $\omega' \in S^{p'-1}$,
 $\omega'' \in S^{p''-1}$,
 $\eta' \in S^{q'-1}$,
 $\eta'' \in S^{q''-1}$,
 $r \ge 0$,
 and
 $\psi \in \Bbb R$.
Then our interest is in a neighbourhood of $r=0$.
In this coordinate,
 we have
$$
   M_+ = \{r \neq 0, \cos 2 \psi > 0\}
$$
The Jacobian matrix of the transform
 $(\theta, \varphi) \to (r, \psi)$
 is given by
$$
     \pmatrix
               \frac {\partial r}{\partial \theta}
             & \frac {\partial r}{\partial \varphi}
\\
               \frac {\partial \psi}{\partial \theta}
             & \frac {\partial \psi}{\partial \varphi}
       \endpmatrix
     =
       \pmatrix
                - \cos \psi \sin \theta
             &  - \sin \psi \sin \varphi
\\
               \frac1{r} {\sin \psi \sin \theta} 
             & \frac{-1}{r} {\cos \psi \sin \varphi}
       \endpmatrix.  
$$
\proclaim{Lemma~\num}
{\rm 1)}\enspace
The standard measure on $S^{p-1} \times S^{q-1}$
 is locally represented as
$$
     \text{smooth function of }(r,\psi,\omega',\eta',\omega'', \eta'')
    \times r^{p' + q' -1} \
     d r d \psi d \omega' d \eta' d \omega'' d \eta''.  
$$
{\rm 2)}\enspace
Any smooth vector field on $S^{p-1} \times S^{q-1}$ near $r = 0$
 is a linear combination of
$$
  \frac{1}{r \cos \psi} X', 
   X'',
  \frac{1}{r \sin \psi} Y', 
  Y'',
  \der{r}, 
  \frac{1}{r}\der{\psi}, 
$$
  whose coefficients are smooth functions of 
 $(\omega',\omega'', \eta',\eta'', r, \psi)$.
Here,
 $X', X'', Y'$ and $Y''$ are smooth vector fields on 
 $S^{p'-1}, S^{p''-1}, S^{q'-1}$ and $S^{q''-1}$,
 respectively.
\endproclaim
  
By noting the relations
$$
\frac {\cos^2 \theta + \cos^2 \varphi}{\cos^2 \theta - \cos^2 \varphi}
     = \frac 1 {\cos 2 \psi},
\quad
     \frac {\sin^2 \varphi + \sin^2 \theta}{\sin^2 \varphi - \sin^2 \theta}
      =\frac {2-r^2}{r^2 \cos 2 \psi},
$$
 $T_+(f_1 f_2) (u,v)$ is locally written as 
$$
  \align
   &B \  (r^2 \cos 2 \psi)^{-\frac {p+q-4} 4}
   (\cos 2 \psi)_+^{\frac {2 \lambda'+ p'+q'-2} 4}
   (r^2 \cos 2 \psi)_+^{\frac {2 \lambda'' + p'' + q''-2} 4}
\\
   =&
   B \ 
   r_+^{\frac {2 \lambda'' - p' -q'+2} 2}
   (\cos 2 \psi)_+^{\frac {\lambda' +\lambda''} 2}, 
  \endalign
$$
 where $B$ is a smooth function of variables
 $\omega',\omega'', \eta',\eta'', r_+, (\cos 2\psi)_+^{\frac12}$.
Therefore,
 by using Lemma~\num,
 we have:
$$
\alignat4
&2\lambda'' > -1, 
&&\lambda' + \lambda'' > -1
 &\quad \Rightarrow &\quad T_+ f \in L^2_{\roman loc},&&
\\
&2 \lambda'' \ge 0, 
&&\lambda' + \lambda'' \ge 1
 &\quad \Rightarrow &\quad 
 Y(T_+ f) 
    \in L^{2-\epsilon}_{\roman loc} 
&& \quad \text{ for any } \ \epsilon > 0, Y \in \frak X(M),
\\
&2 \lambda'' > 2-p'-q', 
&&\lambda' + \lambda'' > 2 
&\quad \Rightarrow &\quad
 Y_1 Y_2(T_+ f)    \in L^{1}_{\roman loc}
&& \quad  \text{ for any } \ Y_1, Y_2 \in \frak X(M),
\endalignat
$$
 in a neighbourhood of the boundary  $\partial M_+$ for Case (2).

The asymptotic estimate for Case (3) parallels to that of Case (2).

\remark{Remark \num}
Assume $\lambda'= \lambda'' \in \Bbb Z + \frac {p'+q'} 2$.
Then $r_+^{\frac {2 \lambda'' - p' -q' +2}{2}}$ is bounded near $r = 0$
 if and only if $\lambda' \ge \frac {p' + q'} 2 -1$,
 equivalently,
  $\lambda' > \frac {p' + q'} 2 -2$,
 which means that
 $\Bbb C_{\lambda'}$ is in the good range
 with respect to the $\theta$-stable parabolic subalgebra
 defined by $\Bbb C_{\lambda'}$ in the sense of Vogan.  
\endremark

\def\sc{8}
\sec{}
We end with some remarks and conjectures, primarily
concerning the precise form of the discrete spectrum
and also the continuous spectrum (where it would be very
interesting to develop the complete Plancherel formula,
given our explicit intertwining operator).

\noindent{1)}\enspace
(multiplicity free property)

Each irreducible component in Theorem~\ch.1
 occurs as multiplicity free.
If $p=2$,
 then $\varpi^{p,q}$ is the direct sum of an irreducible highest weight
 module and a lowest one.
It was proved in \xkmf\ that
 the multiplicity in 
 the full Plancherel formula (both discrete and continuous spectrum)
 is at most one,
 in the branching law of any highest weight module
 of scalar type with respect to any symmetric pair.

\noindent{2)}\enspace
(at most finitely many discrete spectrum, and full discrete spectrum)

If $p', q', p''$ and $q''$ satisfy
$$
  \text{$\min(p', q'') \le 1$ and $\min(q',p'') \le 1$},
\tag \num.1
$$
 then the parameter set in Theorem~\ch.1 is empty, namely,
$$
\Adisc(p', q') \cap \Adisc(q'',p'')
 =  \Adisc(q', p') \cap \Adisc(p'',q'') = \emptyset.
$$
We conjecture that
 there are at most finitely many discrete spectra
 in the branching law $\varpi^{p,q}|_{G'}$ if (\num.1) holds.
We further conjecture that the full discrete spectrum is as
in Theorem 9.1 with $A'(p,q)$ replaced by $A_0(p,q)$
everywhere. 

\noindent{3)}\enspace
(no discrete spectrum)

Furthermore,
 if we exclude the case such as $G' = O(p,q-1) \times O(1)$
 (see \S 7.2),
 namely,
 if
$$
  \text{$\min(p', q'') \le 1$, $\min(q',p'') \le 1$},
  p' + q' > 1, \text{ and } p''+q'' > 1,
\tag \num.2
$$
 then
$$
A_0(p', q') \cap A_0(q'',p'') =  A_0(q', p') \cap A_0(p'',q'') = \emptyset.
$$
It is likely that there is no discrete spectra
 in the branching law $\varpi^{p,q}|_{G'}$ if (\num.2) is satisfied.

We note that the condition (\num.2)
 is equivalent to that
 at least one of $X(p', q')$ or $X(q'', p'')$
 is a non-compact Riemannian symmetric space
 and 
 at least one of $X(q', p')$ or $X(p'', q'')$
 is a non-compact Riemannian symmetric space.

\noindent{4)}\enspace
(discretely decomposable case)

The opposite extremal case is when
$$
 \min(p', p'', q', q'') = 0.
\tag \num.3
$$
As we have proved in Theorem~4.2,
 the restriction $\varpi^{p,q}|_{G'}$ is discretely decomposable
 without any continuous spectrum.
We have obtained the full branching formula in Theorem~7.1
 by using Theorem~4.2 and the $K$-type formula of $\varpi^{p,q}$.

If we employ only the method in this section to the special case (\num.3),
 then we do not have to consider Cases~(2) and (3) in \S 9.7.
Then, Theorem~9.1 exhausts all discrete spectra in Theorem~7.1
 in most cases, but there are some few exceptions.
To be precise,
 we consider the case $p'' = 0$ without loss of generality.
Then,
 in view of Theorem~7.1,
 the right side of Theorem~\ch.1 exhausts all discrete spectra
 if $q'' \ge 5$;
 while
 at most two of $(\frak g', K')$-modules
 $\pip{p}{q'}{\lambda} \boxtimes \pim{0}{q''}{\lambda}$ are missing
 in Theorem~\ch.1
 if $q' \le 4$. %
The precise missing parameters in the case $p'' = 0$ and $q' = 0$ are:
  $\lambda = \pm \frac{1}{2}$  ($q'' = 1$);
  $\lambda = 0, 1$  ($q'' = 2$);
  $\lambda =  \frac{1}{2}$  ($q'' = 3$);
  and $\lambda = 1$  ($q'' = 4$).
In order to 
 to cover all missing parameter
 by the purely geometric method of this section, 
 we should notice a specific feature in the case $p'' = 0$: 
\item{i)} 
 $M_+$ is a dense subspace of $M$.
\item{ii)} 
Any real analytic functions on $M_+$ satisfying the Yamabe equation
 corresponding to $K'$-finite vectors of the $(\frak g', K')$-module
 with the above missing parameter
 extend to real analytic functions on $M$ (see Corollary~4.3~(2)).

\noindent{5)}\enspace
(explicit continuous spectrum)
\def \princeprime#1#2#3#4{{#1}\text{-}\operatorname{Ind}_{\Pmax_{#4}}^{G'_{#4}}
({#3} \otimes \Bbb C_{#2})}

We conjecture that
$$
 \princeprime{L^2}{\sqrt{-1}\lambda}{\epsilon}{1} 
 \boxtimes
 \princeprime{L^2}{\sqrt{-1}\lambda}{\epsilon}{2}
 \quad (\lambda \in \Bbb R)
$$
 is a continuous spectrum with multiplicity free
 if $\min(p', p'', q', q'') > 0$.

\redefine\cite{{}}

\enddocument